\documentclass[11pt,letterpaper]{amsart}
\usepackage[utf8]{inputenc}
\usepackage{enumitem,kantlipsum}
\usepackage[utf8]{inputenc}
\usepackage{amsmath}
\usepackage{marginnote}
\usepackage{amsfonts}
\usepackage{color}
\usepackage{tikz}
\usetikzlibrary{patterns}
\usepackage[color,matrix,arrow,all,cmtip,rotate]{xy}
\usepackage{stmaryrd}
\usepackage{amsthm}
\usepackage{mathrsfs}
\usepackage{amssymb}
\usepackage{hyperref}
\usepackage{graphicx}
\usepackage{wrapfig}
\usepackage{float}

\graphicspath{{dibujos/} }

\numberwithin{equation}{section}

\newtheorem{theorem}{Theorem}[section]
\newtheorem{lemma}[theorem]{Lemma}
\newtheorem{proposition}[theorem]{Proposition}

\newtheorem*{proposition*}{Proposition}
\newtheorem*{theorem*}{Theorem}
\newtheorem*{corollary*}{Corollary}

\newtheorem{corollary}[theorem]{Corollary}
\newtheorem{definition}[theorem]{Definition}
\newtheorem{assumption}{Assumption}

\theoremstyle{remark}
\newtheorem{remark}[theorem]{Remark}
\newtheorem{example}[theorem]{Example}

\newcommand{\tF}{\tilde{F}}
\newcommand{\lF}{\ol{F}}
\newcommand{\tN}{\tilde{N}}
\newcommand{\tvarphi}{\tilde{\varphi}}

\newcommand{\Fix}{{\bf F}_{\ol{n},\ol{d}}}
\newcommand{\HWFix}{\mathbf{W}^{\ol{m}}\Fix}
\newcommand{\HWFixd}{\mathbf{W}^{\ol{m}}_{\ol{d}}\Fix}
\newcommand{\Fixd}{{\bf F}_{\ol{n},\ol{d}}}

\newcommand{\IW}{\mathrm{I}_W^{n_0,n_1}}
\newcommand{\IC}{\mathrm{I}_C^{n_0,n_1}}
\newcommand{\nmin}{n_{min}}
\newcommand{\nmax}{n_{max}}
\newcommand{\supp}{\mathrm{supp}}
\newcommand{\Ker}{\mathrm{Ker}}
\newcommand{\Coker}{\mathrm{Coker}}
\newcommand{\gr}{\mathrm{gr}}
\renewcommand{\Im}{\mathrm{Im}}
\newcommand{\Pic}{\mathrm{Pic}}

\newcommand{\Ext}{\mathrm{Ext}}
\newcommand{\N}{\mathrm{N}}
\newcommand{\M}{\mathrm{M}}
\newcommand{\GL}{\mathrm{GL}}

\newcommand{\U}{\mathrm{U}}
\newcommand{\End}{\mathrm{End}}

\newcommand{\Aut}{\mathrm{Aut}}
\newcommand{\Sym}{\mathrm{Sym}}
\newcommand{\codim}{\textrm{codim}}
\newcommand{\rk}{\textrm{rk}}

\newcommand{\RR}{\mathbb{R}}
\newcommand{\CC}{\mathbb{C}}
\newcommand{\ZZ}{\mathbb{Z}}
\newcommand{\HH}{\mathbb{H}}
\newcommand{\PP}{\mathbb{P}}

\newcommand{\EE}{\mathbb{E}}

\newcommand{\Ss}{\mathcal{S}}
\newcommand{\Ff}{\mathcal{F}}
\newcommand{\Oo}{\mathcal{O}}

\newcommand{\Ee}{\mathcal{E}}

\newcommand{\Hh}{\mathcal{H}}

\newcommand{\Tt}{\mathcal{T}}

\newcommand{\Wd}{\mathbf{W}^{\delta}_{\ol{n}}}
\newcommand{\Zd}{\mathbf{Z}^{\delta}_{\ol{n}}}
\newcommand{\Fd}{\mathbf{F}^{\delta}_{\ol{n}}}
\newcommand{\bfW}{\mathbf{W}}
\newcommand{\bfC}{\mathbf{C}}
\newcommand{\Cd}{\mathbf{C}^{\delta}_{\ol{n}}}

\newcommand{\ol}[1]{\overline{#1}}

\begin{document}
\title[Wobbly moduli of chains]{Wobbly moduli of chains,\\equivariant multiplicities and $\U(n_0, n_1)$-Higgs bundles}
\author{Ana Pe\'on-Nieto}
\address{Ana Pe\'on-Nieto,
	School of Mathematics,  University of Birmingham, Watson Building, Edgebaston, Birmingham B15 2TT, UK}
\email{a.peon-nieto@bham.ac.uk}
\thanks{Project funded through the schemes Proyectos de Consolidación Investigadora (grant nº CNS2022-136042) and Proyectos de generación del conocimiento (nº PID2023-147785NA-I00).}
\date{}

\begin{abstract}
    We give a birational description of the reduced schemes underlying the irreducible components of the nilpotent cone and the $\CC^\times$-fixed point locus of length two in the moduli space of Higgs bundles. Using these results, we prove Drinfeld's conjecture for the sublocus of type $(n_0,n_1)$ fixed points. We introduce the notion of $\U(n_0,n_1)$-wobbliness (stronger than the one of wobbliness) and show that fixed point components of type $(n_0,n_1)$ are wobbly in rank higher than three, if and only if they are also $\U(n_0,n_1)$-wobbly. This yields a computable criterion to check wobbliness of fixed point components, that simplifies the existing ones. We analyse the virtual equivariant multiplicities of fixed points  of type $(n_0,n_1)$ and their Euler pairings with downward flows for type $(1,\dots, 1)$ fixed points. We find that both invariants fail to fully detect all wobbly components for ordered partitions other than $(2,1)$.  
\end{abstract}

--\maketitle
\tableofcontents
\begin{flushright}
{\small{\it To Candela and In\'es}} \\ \vspace*{0.25cm}
\end{flushright}
\section{Introduction}
The moduli space $\M_X(n,d)$ of Higgs bundles of rank $n$ and degree $d$ on a Riemann surface $X$ of genus at least two has been an object of intense study since its introduction by Hitchin  35 years ago \cite{SDE}. In spite of their already long history, they keep proving central through many applications both in geometry and theoretical physics. To cite a few examples, Higgs bundles appear in relation with integrable systems \cite{Duke}, mirror symmetry in its many forms \cite{HT, KW, GWZ, HEnhanced}, quantisation \cite{Gukov, GW}, and the geometric Langlands programme \cite{DPLectures,DPHecke,DPS}. 

In this article, we focus on some key objects in the geometry of moduli spaces of Higgs bundles. These are the fixed points of a natural $\CC^\times$-action which were first studied by Simpson \cite{Simpson}, and later described in terms  of moduli spaces of chains \cite{Chains}. Such points  determine the topology of $\M_X(n,d)$. Indeed, taking limits to zero retracts $\M_X(n,d)$ to the nilpotent cone, and, more precisely, to the fixed points  of the $\CC^\times$-action.
 Other interesting objects that appear from fixed points via the $\CC^\times$ flow include the irreducible components of the nilpotent cone. These can be constructed as upward flows to  fixed point components (that is, taking limits at infinity to a given fixed point component). Another instance, meaningful in mirror symmetry, is the production of $\CC^\times$-invariant branes such as downward flows to (or upward flows from, see Remark \ref{rk:downward_vs_upward}) the so called  very stable Higgs bundles \cite{HH,HEnhanced}.  This notion, very stability, and its opposite, wobbliness, are at the core of many central problems, such as the determination of multiplicities of the components of the nilpotent cone \cite{HH,HitMult}, the definition of divisors inside $\M_X(n,d)$ \cite{Laumon,PalPauly,PalDrinfeld}, and the already cited mirror symmetry of Hitchin systems \cite{HT,HEnhanced} and geometric Langlands via abelianization of Higgs bundles \cite{DPHecke,DPLectures}. This explains the increasing interest in such objects. After the study of fixed point components of nilpotent order $n$ by Hausel and Hitchin, the present paper contributes by analysing the wobbliness of fixed points of nilpotent order two. This is the first step in a programme whose final goal is to determine the components with very stable points, those to which the existing techniques apply, as well as to increase the range of examples. As explained above, our work is motivated by the role played by wobbly and very stable bundles in several important strands of research. 
 
Indeed, the study of very stable bundles can be traced back to the work of Laumon \cite{Laumon} on the nilpotent cone of the moduli space of Higgs bundles. Laumon, following Drinfeld, introduced very stable bundles, and proved that they form a dense open set in the moduli space of stable bundles. He announced the nowadays called Drinfeld conjecture. According to this, the complement of very stable bundles (subsequently named ``wobbly'' by Donagi--Pantev \cite{DPLectures}) is of pure codimension one.  This result  was proven by Pal--Pauly \cite{PalPauly} in rank two, and announced by Pal \cite{PalDrinfeld} for arbitrary rank. The rank three case is studied in \cite{PPrk3}.

Wobbly bundles are also crucial in relation with geometric Langlands from abelianisation of Higgs bundles \cite{DPLectures}. According to Donagi--Pantev, the right setup towards this programme involves the study of parabolic Higgs bundles on the moduli space of vector bundles minus the ``shaky" divisor. They conjectured the equality of the wobbly and shaky loci, a proof of which in the smooth moduli space was provided by the author \cite{W=S}. A toy model of Donagi--Pantev's programme was produced by Hausel--Hitchin \cite{HH} using fixed points of maximal nilpotent order rather than minimal.

More recently, criteria for wobbliness of fixed points in terms of properness of the Hitchin map \cite{PPvs,Hacen, HH} has opened the  way to the computation of multiplicities of the irreducible components of the nilpotent cone, only known until then for the Hitchin section \cite{BR}, in terms of the degree of a restriction of the Hithin map \cite{HH},  computed for vector bundles in \cite{BNR}. Indeed, downward flows to very stable fixed points are proper subvarieties of the moduli space, intersecting the nilpotent cone with generic multiplicity \cite{HH}. Now, although very stability is generic, it can be empty \cite{HH,PPrk3}. Components with no very stable points are called wobbly, and very stable otherwise. The current knowledge of the geometry of very stable components is thus deeper than that of wobbly components. The determination of which are which is a basic problem necessary in order to understand fundamental questions such as the structure of the nilpotent cone and the related dynamics of the $\CC^\times$-flow, amongst others.

An obstruction to very stability of components is provided by Hausel and Hitchin \cite{HH}. They define an invariant, called virtual equivariant multiplicity, which for very stable components recovers the actual multiplicity. These invariants are power series which are polynomial when the component contains a very stable point. Their failure to be such therefore implies wobbliness of the given component. Nonetheless, this failure does not characterise wobbliness, as virtual equivariant multiplicities are known to be polynomial also for some concrete wobbly components (e.g., of type $(3,1)$ \cite{HH}).   

In this paper,  we show that all fixed point components of nilpotent order two are wobbly in rank higher than 3, in contrast to what happens in the maximal order case \cite{HH}. We also compute two invariants: the virtual equivariant multiplicities of type $(n_0,n_1)$ and a symmetric version of the equivariant Euler class pairing the $BAA$ branes given by downward flows to type $(n_0,n_1)$ and $(1,\dots,1)$ fixed points. We analyse  the detection power for wobbliness of both invariants.  We find that, for all ranks except $3$, and all partitions $n=n_0+n_1$ of length two, there exist wobbly components of type $(n_0,n_1)$ for which the virtual equivariant multiplicities are polynomial.  These results explain several observations in low rank \cite{HH, PPrk3}. 

Another line we pursue is the study of these questions from the point of view of real groups. In his foundational work, Simpson \cite{Simpson} proved that fixed points are in fact Higgs bundles for real forms. This has been exploited in both directions: results about moduli spaces of real forms have been deduced from those for fixed points and viceversa. For example, the determination of  irreducible components of the fixed point locus allowed to prove connectedness of $\U(p,q)$-Higgs bundle moduli spaces by Bradlow--Garc\'ia-Prada--Gothen--Heinloth \cite{IrredChains}. Very recently, real forms have been used to produce bounded invariants for moduli spaces of chains by Biquard--Collier--Garc\'ia-Prada--Toledo \cite{Arakelov-Milnor}. This interplay between fixed points and real forms points to the basic question of determining to which extent phenomena such as wobbliness of fixed point components are controlled by the associated real form. In this paper we introduce the notion of $G_\RR$-very stability and prove that very stability of length two components is equivalent to the a priori weaker notion of $\U(n_0,n_1)$-very stability. Related notions have been introduced and analysed by Gonz\'alez--Hausel \cite{GH}.

The structure of the present paper is as follows. Section \ref{sec:prelim} contains well known results, with the exception of Proposition \ref{prop:vs_are_dense}, that proves density or emptiness of very stable bundles in a given non strictly semistable fixed point component. In Section \ref{sec:birat}, a birational description of the reduced schemes underlying the fixed points and the associated components of the nilpotent cone with generic Higgs field of nilpotent order two (Theorem \ref{thm:description_fixed_points}).  The fixed points therein have underlying vector bundle of the form $F_0\oplus F_1$, where $\rk(F_i)=n_i$ and Higgs field $\varphi\in H^0(F_1^*F_0K)$. Our proof uses the Brill--Noether theoretic results by Russo--Teixidor i Bigas \cite{RT}. In particular, we determine non emptiness of the fixed point components in terms of an invariant $\delta$ directly deducible from the Toledo invariant from the theory of  $\U(n_0,n_1)$-Higgs bundles by Bradlow--Garc\'ia-Prada--Gothen \cite{BGGP} (albeit more natural for Brill-Noether theoretic reasons). This gives an alternative construction to the one in \cite{triples} in terms of moduli of triples.  As an application, we recover minimal dimensionality of the space of sections $H^0(F_1^*F_0K)$ in Corollary \ref{cor:BN} (already known for $n_0>n_1$ through Brill--Noether theory). All these results yield the necessary conditions for components of the wobbly locus of codimension one to exist. In Section \ref{sec:wobbly_divisors}, we show Drinfeld's conjecture on pure codimensionality of the wobbly locus for a sublocus of wobbly locus, namely, those vector bundles admitting a nilpotent Higgs field of order two (Theorem \ref{thm:wobbly_divisors}). Although we are able to prove non emptiness of this wobbly divisor and its pure codimensionality, we fail to identify the set of invariants linking the irreducible components of the nilpotent cone and those of the wobbly divisors.   

In Section \ref{sec:criteria},  we give some preliminary results towards the next section, as well as criteria for wobbliness in some simple cases. More precisely, we focus on studying the existence of nilpotent order two Higgs fields in downward flows. The main result is contained in Section \ref{sec:wobbly_comps}, where in Theorem \ref{thm:wobbly_comps} we prove wobbliness of fixed point components, and explain some observations of Hausel--Hitchin for, e.g., ordered partitions of type $(3,1)$. After introducing the notion of $G_\RR$-very stability (Definition \ref{def:real_vs}), we prove the equivalence between wobbliness of length two components and the a priori stronger notion of $\U(n_0,n_1)$-wobbliness (Corollary \ref{cor:wobblyiffrealwobbly}). This simplifies the identification of wobbly components by restricting the flows to be considered. In Section \ref{sec:multiplicities} we move on to the computation of the virtual equivariant multiplicities of fixed points (Proposition \ref{prop:equivar_mult}) and determine the range for which wobbliness of the components is captured by the invariants, i.e., the range for which the invariants are not polynomial (Corollary \ref{thm:mults}). In a final Section \ref{sec:Euler_pairings} we analyse invariants related to the equivariant Euler pairings, Euler characteristics of brane homomorphisms. After a general computation of these invariants (Proposition \ref{prop:equiv_Euler_pairings}), we check in Example \ref{ex:euler} that, just as virtual equivariant multiplicities, these new invariants fail to fully detect wobbliness. 
\\\\
\textbf{Acknowledgements.} Many thanks to P. Gothen, T. Hausel, and C. Pauly, for meaningful remarks and discussions. Thanks also to M. Teixidor i Bigas for discussing some aspects of Brill--Noether theory, and M. Gonz\'alez and M. Mazzocco for their comments on a preliminary version of this manuscript. 

\section{Preliminaries and notation}\label{sec:prelim}
\subsection{Higgs bundles, the Hitchin map, and fixed points}

Consider the moduli space $\M_X(n,d)$ of semistable Higgs bundles  of rank $n$ and degree $d$ on a Riemann surface $X$ of genus at least two.  Its closed points are $S$-equivalence classes of pairs $(E,\varphi)$ where $E$ is a vector bundle on $X$ and $\varphi\in H^0(\End(E)\otimes K)$ where $K$ is the canonical bundle of $X$. Recall the definition of semistability:
\begin{definition}
The slope of 
a vector bundle $E$ of degree $\deg(E)$ and rank $\rk(E)$ is the quotient \[\mu(E)=\frac{\deg(E)}{\rk(E)}.\]
A Higgs bundle $(E,\varphi)$ is stable (resp. semistable) if for every subbundle $0\subsetneq F\subsetneq E$ such that $\varphi(F)\subseteq F\otimes K$, it holds that
\[
\mu(F)< \mu(E)\,\qquad  (\mathrm{resp.} \mu(F)\leq \mu(E)).
\]
\end{definition}
We will denote S-equivalence of two (Higgs) bundles $E_1,\, E_2$ (resp. $(E_1,\varphi_1),\ (E_2,\varphi_2)$) by $E_1,\sim_S E_2$ (resp. $(E_1,\varphi_1)\sim_S (E_2,\varphi_2)$). We will abuse notation by identifing a semistable Higgs bundle $(E,\varphi)$ with the point it defines in the moduli space. We will denote by $\gr(E,\varphi)$ the unique polystable representative in its S-equivalence class (as it corresponds with the graded Higgs bundle for the Jordan--H\"older filtration).

Let 
\begin{equation}
    \label{eq:HitchinMap}
h:\M_X(n,d)\longrightarrow B_n=\bigoplus_{i=1}^n H^0(K^i)\qquad (E,\varphi)\mapsto \det(xId-\varphi)
\end{equation}
be the Hitchin map. It is a proper map whose fibers are complex Lagrangians of pure dimension. The fiber over zero, the so called global nilpotent cone $h^{-1}(0)$, consists of Higgs bundles with nilpotent Higgs fields.

There is a natural $\CC^\times$-action on the moduli space of Higgs bundles 
\[
t\cdot(E,\varphi)=(E,t\varphi).
\]
The Hitchin map is $\CC^\times$-equivariant for a suitable weighted action on $B_n$. This, together with properness of the Hitchin map, implies the existence of a limit 
\[\lim_{t\longrightarrow 0}(E,\varphi)\in h^{-1}(0),
\]
which is moreover a fixed point for the $\CC^\times$-action. Similarly, equivariance of the Hitchin map and properness of the fibers imply that for any $(E,\varphi)\in h^{-1}(0)$
\[
\lim_{t\longrightarrow \infty} (E,\varphi)\in h^{-1}(0)^{\CC^\times},
\]
where $ h^{-1}(0)^{\CC^\times}$ denotes the fixed point set of the $\CC^\times$-action (necessarily nilpotent). These were classified by Simpson \cite{Simpson}. They are of the form
\[
\Ee=\left(\bigoplus_{i=0}^s F_i,\bigoplus_{i=i}^s \phi_i\right)
\]
where $F_i$ is a rank $n_i$ degree $d_i$ vector bundle (where $\ol{n}=(n_0,\dots, n_s)$ and $\ol{d}=(d_0,\dots, d_s)$ are ordered partitions of $n$ and $d$ respectively), and $\phi_i\in H^0(X,F_{i+1}^*F_iK)$. 
\begin{definition}\label{def:type_fixed_point}
The Simpson type (type for short, see Remark \ref{rk:Laumon_vs_Simpson_type}) of the fixed point $\Ee=(\bigoplus_{i=0}^s F_i,\bigoplus_{i=i}^s \phi_i)$ is the pair $(\ol{n},\ol{d})$, where $\ol{n}$ and $\ol{d}$ are ordered partitions of $n$ and $d$ respectively of the same length, corresponding to the ranks and degrees of the graded terms $F_i$. We will denote by $\mathbf{F}_{\ol{n},\ol{d}}$ an irreducible component of the fixed point set of given type $({\ol{n},\ol{d}})$.
\end{definition}
\begin{remark}\label{rk:Laumon_vs_Simpson_type}
In \cite[Def. 1.7]{Laumon}, Laumon introduced a different notion of (nilpotent) type for an arbitrary nilpotent Higgs bundle, henceforth called Laumon-type. Given ordered partitions $(\ol{n},\ol{d})$ with $n_0\geq \dots\geq n_s$ and such that the generic point in the component $\mathbf{F}_{\ol{n},\ol{d}}$ satisfies that the summands $F_i$ are stable, then the Simpson type of a fixed point as per Definition \ref{def:type_fixed_point} is a {\it nilpotent Laumon-type} (thus also a Laumon-type) in Laumon's terminology. 
\end{remark}
 Types determine the irreducible components of the fixed point set \cite[Corollary 3.3]{Bozec}. Moreover, there is a one-to-one correspondence between the irreducible components of the nilpotent cone and of the fixed point set. The relationship is given by the Bialynicki-Birula stratification: taking limits at $\infty$ defines a Zariski locally trivial fibration 
 \[
 {\Fixd}^-\longrightarrow\Fixd
 \]
 where $\Fixd$ denotes an irreducible component of the fixed point locus and \[\Fixd^-:=\{(F,\psi)\in\M_X(n,d)\,:\,\lim_{t\to\infty}t\cdot(F,\psi)\in\mathbf{F}_{\ol{n},\ol{d}}\}.
 \]
 Then, components of the nilpotent cone are precisely the closures  $\mathbf{C}_{\ol{n},\ol{d}}:=\ol{\mathbf{F}^-}_{\ol{n},\ol{d}}\subset h^{-1}(0)$. 
\subsection{Wobbly and very stable points and components}

\begin{definition}\label{def:very_stable}
Let $\Ee$ be a fixed point. Let
\[
\Ee^+=\{(F,\psi)\in\M_X(n,d)\,:\,\lim_{t\to 0}(F,\psi)=\Ee\}.
\]
A stable fixed point is {\bf very stable} if $\Ee^+\cap h^{-1}(0)=\{\Ee\}$. Otherwise, it is called {\bf wobbly}. 
\end{definition}
\begin{remark}\label{rk:downward_vs_upward}
    Hausel and Hitchin refer to $\Ee^+$ as the upward flow from $\Ee$. We will instead use the equivalent term downward flow to $\Ee$.
\end{remark}
This notion of very stability generalises the one by Drinfeld \cite{Laumon} for vector bundles in a natural way, but it differs in a fundamental aspect: while Drinfeld's definition for vector bundles only considers fixed points, Definition \ref{def:very_stable} puts the emphasis in downward flows from the nilpotent cone, which implies that the criterion for very stability via the Hitchin map \cite{PPvs,Hacen,HH} holds. While both are equivalent, the following lemma gives an alternative proof to the one in 
\cite{HH}.
\begin{lemma}\label{lm:wobbly_flow_from_nilp}
    Let $\Ee\in\M_X(n,d)^{\CC^\times}$ be stable. Then, $\Ee$ is very stable if and only if \[\ol{\Ee^+}\cap h^{-1}(0)\setminus h^{-1}(0)^{\CC^\times}\neq \emptyset.\]
\end{lemma}
\begin{proof}
The necessity of the condition is clear, as
\[
\Ee^+\setminus\{\Ee\}\subset \ol{\Ee^+}\cap h^{-1}(0)\setminus h^{-1}(0)^{\CC^\times}\]

For the converse, a slight modification of the proof in \cite{PPvs} works in this setup. We sketch the proof for the reader's convenience. By wobbliness, there is $\Ff\in\ol{\Ee^+}\cap h^{-1}(0)\setminus \Ee$. Let $\iota: C\longrightarrow \ol{\Ee^+}$ be a smooth curve with every point in $\Ee^+$ except for $o\in C$, with $\iota(c)=\Ff$ (which exists by, e.g., \cite[Lemma 2.4]{PPvs}). Then, the $\CC^\times$ action defines a surface  rationally embedded in $\ol{\Ee^+}$ (cf. \cite[Prop. 3.1]{PPvs}) by
\[
\xymatrix{\CC\times C\ar[r]^j&\ol{\Ee^+}\\
(t,p)\ar@{|->}[r]&t\cdot\iota(p).
}
\]
This is well defined beyond $(0,o)$. By the discussion following \cite[Question E]{Hironaka}, a finite number of blowups resolves the morphism (see the discussion following  \cite[Prop. 3.1]{PPvs}), yielding
\[\xymatrix{\widehat{\CC\times C}\ar[r]^j&\ol{\Ee^+}}.
\]
where $\widehat{\CC\times C}$ denotes the resulting scheme after a finite number of blowups. The exceptional divisor  $D=\bigcup \PP^1\subset \widehat{\CC\times C}$ is connected, joining  $\Ee=\lim_{\lambda\to 0}j(\lambda, c)$ and $\Ff=\lim_{c\to o}j(0,c)$. One of the lines in $D$ thus intersects $\Ee^+$, which is dense and open inside $\ol{\Ee^+}$, and thus by irreducibility, the line must contain a dense open subset of elements in $\Ee^+$. These elements are in the nilpotent cone. Indeed, by Hartog's theorem, the induced rational Hitchin map $h:\CC\times C\dasharrow B_n$ extends uniquely to the indeterminacy locus. By continuity, the image is zero, and so $h(D\cap\M_X)=0$ (cf. the discussion following \cite[Prop. 3.1]{PPvs}). On the other hand, at least one of the lines intersecting $\Ee^+$ must be non constant, as otherwise the divisor would not connect $\Ee$ and $\Ff$. This proves the statement.
\end{proof}
It follows from \cite{Laumon,HH}, that very stable fixed points of type $((n),d)$ and type $((1,\dots, 1),(d_1,\dots,d_n))$ are either dense or empty in their corresponding components of the nilpotent cone.  Proposition \ref{prop:vs_are_dense} proves this is always the the case for non singular components.
\begin{proposition}\label{prop:vs_are_dense}
    Let $\emptyset\neq \Fix^s\subset\Fix$ be the stable locus. Then either $\Fix^s$ is fully wobbly or it contains a dense open set of very stable points.
\end{proposition}
\begin{proof}
Consider the downward flow map $\pi: (\Fix^s)^+\longrightarrow \Fix^s$, which is closed and locally of finite type as it is a Zariski locally trivial fibration \cite[Theorem 4.1]{BB} by half dimensional affine spaces \cite[Prop. 2.10]{HH}. Consider the restriction: 
\[
\pi_0: \Fix^+\cap h^{-1}(0)\longrightarrow \Fix\]
which is closed and locally of finite type as it is the composition of two such morphisms. Then, by \cite[Theorem IV.13.1.5]{EGA}, the dimension of the fibers is upper semicontinuous, and thus  generically minimal. By definition, the minimal dimension is $0$ if there is a stable very stable bundle. So, either generic smooth points are very stable, or the component is wobbly along the smooth locus. 
\end{proof}
Proposition \ref{prop:vs_are_dense} justifies the following definition. 
\begin{definition}
A non-singular fixed point component is called very stable if it contains a very stable point. Otherwise, it is called wobbly. 
\end{definition}
Given a $\CC^\times$-module $V$, for $\lambda\in\ZZ$, let $V_\lambda$ denote the $\lambda$ weight space, and let $V^+=\bigoplus_{\lambda>0}V_\lambda$.  Let $\chi_{\CC^\times}$ denote the character of a $\CC^\times$-module. For a fixed point $\Ee$, $T_\Ee\M_X(n,d)$ carries a $\CC^\times$-action, as so does $\Sym({T^*_\Ee})$. Let $T_\Ee^+$ denote the $\CC^\times$ submodule of $T_\Ee\M_X(n,d)$ corresponding to positive weights.

\begin{definition}[\cite{HH}]
The virtual equivariant multiplicity of $\Ee$ is the fraction
\[
m_\Ee(t)=\frac{\chi_{\CC^\times}(\Sym({T^*_\Ee}^+))}{\chi_{\CC^\times}(\Sym(B_n^*))}
\]
where $B_n$ is the Hitchin base, cf. \eqref{eq:HitchinMap}.
\end{definition}
\begin{remark}
When $\Ee$ is very stable, this is a polynomial whose value $m_\Ee(1)$ at one matches the multiplicity of the irreducible component containing the point \cite[Theorem 5.2]{HH}. 
\end{remark}
\subsection{$\U(p,q)$-Higgs bundles}\label{sec:Upq}
Let $G_{\RR}<\GL(n,\CC)$ be a real form. Let $H_{\RR}<G_{\RR}$ be a maximal compact subgroup. Denote by $\mathfrak{g}_{\RR}:=Lie(G_{\RR})$, $\mathfrak{h}_{\RR}:=Lie(H_{\RR})$. Let 
\[
\mathfrak{g}_{\RR}=\mathfrak{h}_{\RR}\oplus \mathfrak{m}_{\RR}
\]
be the Cartan decomposition. Let $H:=H_{\RR}^\CC$, $\mathfrak{h}:=\mathfrak{h}_{\RR}^\CC$,  $\mathfrak{m}:=\mathfrak{m}_{\RR}^\CC$ be the complexifications. 
\begin{definition}
A $G_{\RR}$-Higgs bundle is a Higgs bundle $(E,\varphi)$ where $E=E_H(\CC^n)$ is the associated bundle (for the standard representation) for a principal $H$-bundle $E_H$ and $\varphi\in H^0(E_H(\mathfrak{m})\otimes K)$, where $H$ acts on $\mathfrak{m}$ via the isotropy representation and we identify $E_H(\mathfrak{m})\subset\End(E_H(\CC^n))$ via $\mathfrak{m}\subset\mathfrak{gl}_n(\CC)\cong\End(\CC^\times)$.
\end{definition}
A stability condition can be defined yielding a moduli space $\M_X(G_\RR)$ admitting a natural map
\begin{equation}
    \label{eq:MR_MC}
\M_X(G_\RR)\longrightarrow \M_X(n,d).
\end{equation}
By abuse of notation, we will identify the image of this map with $\M_X(G_\RR)$, althought the fibers of the map \eqref{eq:MR_MC} are non trivial \cite{Abelian}.
\begin{definition} A $\U(p,q)$-Higgs bundle is a Higgs bundle $(E,\varphi)$ where $E=V\oplus W$ where $\rk(V)=p$, $\rk(W)=q$ and $\varphi=(\beta,\gamma)$ with $\beta:V\longrightarrow WK$, $\gamma:W\longrightarrow VK$.
\end{definition}
The ranks $(p,q)$ and degrees $(d_V,d_W)$ of $V$ and $W$ combine in the so-called Toledo invariant, defined by
\begin{equation}\label{eq:Toledo}
\tau=2\frac{qd_V-pd_W}{p+q}
\end{equation}
This invariant is bounded by
\begin{equation}\label{eq:bounds_Toledo}
    0\leq|\tau|\leq 2\min\{p,q\}(g-1).
\end{equation}
The irreducible components of $\M_X(\U(p,q))$ are classified by the degrees $(d_V,d_W)$ except in the maximal Toledo case $\tau=\pm 2\min\{p,q\}(g-1)$, for which the Higgs bundle is generically strictly semistable, with $\mu(F_1\oplus F_1K^*)=\mu(F_0\oplus F_1/F_1\oplus F_1K^*)$. Thus, the S-equivalence class is $(F_1K^*\oplus S\oplus F_1,\varphi)$ for some vector bundles $F_0,\, S,\ F_1$ where $F_0=S\oplus F_1K^*$ and the Higgs field is the composition $\varphi: F_1\stackrel{1}{\longrightarrow} F_1K^*K\hookrightarrow F_0K$ \cite{BGGP}.

Note that $\mathbf{F}_{(n_0,n_1),(d_0,d_1)}\subset\M_X(\U(n_0,n_1))$. Since inside $\M_X(n,d)$
\[\M_X(\U(n_0,n_1))\cong \M_X(\U(n_1,n_0)),\]
to avoid ambiguity with our notation for fixed point components, for given $(n_0,n_1)$, define
\[\nmin:=\min\{n_0,n_1\},\qquad\nmax:=\max\{n_0,n_1\}.\]
Then 
\[\bigcup_{\tiny\bfC_{\ol{d}}\neq \emptyset}\mathbf{F}_{(n_0,n_1),(d_0,d_1)}\cup \bigcup_{\tiny\bfC_{\ol{d}'}\neq \emptyset}\mathbf{F}_{(n_1,n_0),(d_1',d_0')}\subset\M_X(\U(\nmax,\nmin)).\]
Thus $\mathbf{F}_{(n_1,n_0),(d_0,d_1)}\subset\M(\U(p(n_0,n_1),q(n_0,n_1)))$. 

\section{Fixed points of length two and $\U(n_0,n_1)$-Higgs bundles}\label{sec:birat}
Let $\ol{n}=(n_0,n_1)$,  $\ol{d}=(d_0,d_1)$. In this section,  we will consider irreducible components of fixed points $\Fix\subset\M_X(n,d)^{\CC^\times}$. Namely, $\Ee\in \Fix$ is of the form
\[
\Ee=(F_0\oplus F_1,\varphi),\qquad \varphi\in H^0(F_1^*F_0K).
\]
Observe that the map
\[
t:\M_X(n,d)\longrightarrow \M_X(n,-d)\qquad (E,\varphi)\mapsto (E^*,{^t}\varphi)
\]
maps $\Fix$ to $\mathbf{F}_{(n_1,n_0),-(d_1,d_0)}$. So it is enough to study the fixed points for which $n_0\geq n_1$. These have the advantage of generically having $\Ker(\varphi)=F_0$ (see Theorem \ref{thm:description_fixed_points}).

Let us begin by giving a birational description of the reduced scheme 
underlying the irreducible components.

Define the invariant:
\begin{equation}
    \label{eq:def_delta} \delta=\delta(\ol{n},\ol{d}):=\deg(F_1^*F_0K)=d_0n_1-d_1n_0+2n_0n_1(g-1).
\end{equation} 
We will henceforth label fixed point and nilpotent cone components with this invariant (namely $\Fd,\, \Cd$).

Let $\N_X(n_i,d_i)$ denote the moduli space of semistable vector bundles of rank $n_i$ degree $d_i$, and let $\N^s_X(n_i,d_i)$ be the stable locus. Consider  the closed subscheme of $\N_X(n_0,d_0)\times \N_X(n_1,d_1)$ given by
\begin{equation}
    \label{eq:def_Znd'}
{\Zd}'={\left\{
  (F_0,F_1)\in\N_X(n_0,d_0)\times \N_X(n_1,d_1)\,:\,h^0(\mathrm{gr}(F_0)\mathrm{gr}(F_1)^*K)>0\}
\right\}}.\end{equation}
where $\mathrm{gr}(F_i)$ is the graded object associated to the Jordan--H\"older filtration (so that in particular, for $F_i\in\N_X^s(n_i,d_i)$, $\mathrm{gr}(F_i)=F_i$).
This is a family of twisted Brill-Noether loci over $\N_X(n_0,d_0)$, where we recall that the twisted Brill--Noether locus is defined by
\begin{equation} 
\label{eq:twisted_BN}
        BN^0_{m,d}(E):=\{R\in \N^s_X(m,d)\,:\, h^0(R^*E)\neq 0\}.
        \end{equation} 
For $n_0>n_1$, twisted Brill--Noether loci were studied by Russo--Teixidor  \cite{RT}, and some of their intermediate results on extensions of bundles also contemplate the case $n_0=n_1$. Using the latter, we will give a Brill--Noether construction of the irreducible components of the nilpotent cone. We will need to distinguish the cases $n_0=n_1$ and $n_0>n_1$ at parts, as the tools differ slightly.

Let 
\begin{equation}
    \label{eq:IC}
    \IC:=\left\{\delta\left|\begin{array}{l}
         \nmin(\nmax-\nmin)(g-1)\leq\delta< 2n_0n_1(g-1) \\
         \delta\equiv -n_0d+2n_0n_1(g-1)\mod{n}
    \end{array}\right.\right\}
\end{equation}
with strict lower bound if $n_0\neq n_1$.
Likewise, we define 
\begin{equation}
    \label{eq:IW}
    \IW:=\{\delta\in\IC\colon \delta\leq n_0n_1(g-1)+1.\}
\end{equation}
\begin{lemma}\label{lm:toledo_delta}
    Let $F_0,\ F_1$ be vector bundles of ranks $\rk(F_i)=n_i$, and degrees  $\deg(F_i)=d_i$. Then if  $\varphi\in H^0(F_1^*F_0K)$ satifies that $(F_0\oplus F_1,\varphi)$ is a semistable Higgs bundle, then
    \[\delta \in\IC\cup\{ \nmin(\nmax-\nmin)(g-1),2n_0n_1(g-1)\}.\]
Thus, when  $\delta=2n_0n_1(g-1)$, $(F_0\oplus F_1,\varphi)\sim_S(F_0\oplus F_1,0)$. When  $\delta=n_{\min}(n_{\max}-n_{\min})(g-1)$), if $n_0\geq n_1$ (resp. $n_0\leq n_1$), then $(F_0\oplus F_1,\varphi)\sim_S(F_1K^*\oplus S\oplus F_1,\varphi)$ (resp.  $(F_0\oplus F_1,\varphi)\sim_S(F_0K\oplus S\oplus F_0,\varphi)$).
\end{lemma}
\begin{proof}
\begin{remark}
    \label{rk:toledo_vs_delta}
    Note that with our chosen definition
  \begin{equation}
    \label{eq:delta_vs_toledo} \tau=\left\{\begin{array}{ll}
\frac{2}{n}(\delta-2n_0n_1(n_0-n_1))(g-1)&\textrm{ when } n_0\geq n_1\\
-\frac{2}{n}(\delta-2n_0n_1(n_0-n_1))(g-1)&\textrm{ when } n_0\leq n_1.
\end{array}\right.
\end{equation}
By semistability, negative values of $\tau$ correspond to $n_0\geq n_1$ while positive values of $\tau$ correspond to $n_0<n_1$
so that maximality of the absolute value of $\tau$ implies minimality of $\delta$ when $n_0\geq n_1$ (and viceversa when $n_0\leq n_1$) and minimality of the absolute value of $\tau$ implies maximality of $\delta$ when $n_0\geq n_1$ (and viceversa when $n_0\leq n_1$) in the given ranges. 
\end{remark}
Some of the following arguments can be found in \cite[Theorem 3.32]{BGGP}. 

Semistability of the Higgs bundle implies $\mu(F_0)\leq \mu(F_1)$, or equivalently $\delta\leq 2n_0n_1(g-1)$.  When $\delta=2n_0n_1(g-1)$, then $\mu(F_0)=\mu(F_1)$, so that $(F_0\oplus F_1,\varphi)\sim_S(F_0\oplus F_1,0)$ for any $\varphi\in H^0(F_1^*F_0K)$. Assume $n_0\geq n_1$. Then, if $\delta< n_1(n_0-n_1)(g-1)$, then $\mu(E)<\mu(F_1K^*\oplus F_1)$, violating semistability. Equality $\delta=n_1(n_0-n_1)(g-1)$ implies $\mu(E)=\mu(F_1K^*\oplus F_1)$, so it follows that  $(F_0\oplus F_1,\varphi)\sim_S(F_1\oplus F_1K^*\oplus S,\varphi)$ for some vector bundle $S\cong F_0\oplus F_1/F_1K^*\oplus F_1$. In particular, $\varphi$ has maximal rank. 

Similar arguments yield the same range for $n_0\leq n_1$. Alternatively, note that in that case $F_0\oplus F_1$ has invariant 
\[\delta(\ol{n},\ol{d})=-n_1(-d_0)+n_0(-d_1)+2n_0n_1(g-1)=\delta((n_1,n_0),(-d_1,-d_0))\]
which is the invariant of the dual Higgs bundle. 
\end{proof}

In what follows, we will assume $n_0\geq n_1$. The case $n_1<n_0$ follows by dualisation from the case $n_0>n_1$ (cf. Corollary \ref{cor:comps_n_0<n_1}).

\begin{lemma}\label{lm:sheaf_n0=n1}
Let $n_0=n_1$. There exists a rational projective bundle
\[
\xymatrix{\PP(\tilde{\mathcal{H}}_\delta)\ar@{-->}[r] &\Sym^{\frac{\delta}{n_0}}X\times \N_X(n_0,d_0)}
\]
with fiber $\PP(H^0(F_0^*\otimes\Oo_D))$ over $(D,F_0)\in \Sym^{\frac{\delta}{n_0}}X\times \N_X(n_0,d_0)$. This bundle is the projectivization of a vector bundle with fiber $H^0(F_0^*\otimes\Oo_D)$.
\end{lemma}
\begin{proof}
Consider the projective universal bundle
$\PP(\mathcal{F}_0)\longrightarrow X\times\N_X^s(n_0,d_0)$. Note that \'etale locally there exists a bundle $\mathcal{F}_0$ whose projectivisation is the \'etale local restriction of $\PP(\mathcal{F}_0)$ (this follows from Mumford's general machinery \cite{Mumford} and Luna's slice theorem \cite{Luna}, see \cite[Thm 1.21--1.24]{SimpsonModuli} for a detailed discussion). Now, there is an incidence divisor $W\subset \Sym^n(X)\times X$ such that $W|_{D\times X}\cong D$ as a subscheme of $X$. 
Consider the rational sheaf 
$\mathcal{Y}=\PP(\Oo_{W}\boxtimes \mathcal{F}_0^*)$, which is well defined over the stable locus. Now,  let  $\{X\times U\}\twoheadrightarrow X\times\N_X(n_0,d_0)$ be an \'etale open covering over which the universal family exists. Denote by $p_1$, $p_2$, $p_3$ the projections from $\Sym^{\frac{\delta}{n_0}}X\times X\times U$ to $\Sym^{\frac{\delta}{n_0}}X$, $X$ and $U$ respectively, and define the sheaves 
$\mathcal{Y}_U:=\Oo_{W}\boxtimes \mathcal{F}_0^*|_U:=(p_1\times p_2)^*\Oo_W\otimes (p_1\times p_3)^*\mathcal{F}_0^*|_U$ on 
\[
 Y_U:=\Sym^{\frac{\delta}{n_0}}X\times (X\times U).
\]
Let
 \[
 R_U:=\Sym^{\frac{\delta}{n_0}}X\times (X\times U)\times H^0(Y_U,\mathcal{Y}_U)
\]
\[
 \ol{R}_U:=\Sym^{\frac{\delta}{n_0}}X\times (X\times U)\times \PP(H^0(Y_U,\mathcal{Y}_U)\oplus \CC)
\]
and let $R_U\stackrel{i}{\hookrightarrow}\ol{R}_U$ be the open immersion. 
Let $q_1$, $q_2$, denote the projections to $\Sym^{\frac{\delta}{n_0}}X\times X\times U$ and $\PP(H^0(Y_U,\mathcal{Y}_U)\oplus \CC)$ respectively. Abusing notation, we will also denote by $\mathcal{Y}_U$ the pullback $q_1^*\mathcal{Y}_U$. 

Consider the coherent sheaf $\Oo_{H^0,U}':=i_*(q_2\circ i)^*\Oo_{H^0(\mathcal{Y}_U)}$, whose global sections are precisely $H^0(\mathcal{Y}_U)$. Then 
\[
\Oo_{H^0,U}'\longrightarrow \ol{R}_U
\]
satisfies that the restriction to $X\times \{F_0,D, f\}$ is precisely the section
\[
 f:F_0\longrightarrow \Oo_D.
\]
Now, let $\pi_U:\ol{R}_U\longrightarrow \Sym^{\frac{\delta}{n_0}}X\times U$ denote the projection. Then, by \cite[Theorem III.3.2.1]{EGA},  the sheaf $\tilde{\Hh}_{\delta,U}:=R^0\pi_* \Oo_{H^0}'$ is coherent,  with fiber over $(F_0,D)$ equal to $H^0(F_0^*\otimes \Oo_D)$. Since $h^0(F_0^*\otimes \Oo_D)$ is generically constant over $U$, it follows that $\{\PP(\tilde{\Hh}_{\delta,U})\}$ defines a rational projective bundle satisfying the assumptions. 
\end{proof}
\begin{corollary}\label{cor:compZnn}
    For $n_0=n_1$, ${\Zd}'$ contains an irreducible component $\Zd$ given by  the closure of the image of $\PP(\Hh_\delta)$ under the rational map 
    $r:(F_0,D,\CC^\times f)\mapsto (F_0,\Ker(\CC^\times f)K)$. 
\end{corollary}

\begin{lemma}\label{lm:dimZ}
Let $\Zd\subset {\Zd}'$ be the subscheme defined in Corollary \ref{cor:compZnn} if $n_0=n_1$, and the scheme defined in \eqref{eq:def_Znd'} for $n_0>n_1$. 

Then,  whenever $\delta\in\IC$, $\Zd$ is non empty, irreducible of dimension 
    \begin{equation}
        \label{eq:dimZ}
        \dim \Zd=\left\{\begin{array}{ll}
      (n_0^2+n_1^2)(g-1)+2   &\textrm{ if } \delta\in\IC\setminus\IW\\
         (n_0^2+n_1^2-n_0n_1)(g-1)+\delta+1  &\textrm{ if } \delta\in\IW.
    \end{array}\right.
    \end{equation}
    
    If $\delta>n_0n_1(g-1)$, then $Z'=Z=\N_X(n_0,d_0)\times \N_X(n_1,d_1)$. 
\end{lemma}
\begin{proof}
Non emptiness for $\delta\geq n_1(n_0-n_1)(g-1)$ follows from \cite[Theorem 0.2]{RT}. 
For $n_0=n_1$, using Corollary \ref{cor:compZnn}, the dimension is computed in the proof of \cite[Theorem 0.3]{RT}. 

For $n_0>n_1$, the statement can be deduced from the proof of \cite[Prop. 1.4]{RT}. Indeed, following the arguments in {\it loc. cit.}, the universal family of extensions $\mathcal{E}xt^1\longrightarrow\N_X^s(n_1,d_1)\times \N_X^s(n_0-n_1,d_0-d_1+2n_1(g-1))$ is birational to a vector bundle over the dense open set of pairs $(F_1,S)$ such that $H^0(S^*F_1K^*)=0$. This bundle is irreducible, and admits a rational projection $\pi: \PP(\mathcal{E}xt^1)\dasharrow\N_X(n_1,d_1)\times \N_X(n_0,d_0)$ mapping a triple $(F_1,S,f_0)$ to $(F_1,F_0)$. The image of $\pi$ is irreducible, contained in $\Zd$. Moreover, it contains all the pairs in $\Zd$ such that $F_1$ and $S$ are stable. So all we need to check is that  $\Zd$ contains no irreducible component characterised by $F_1$ or $S$ being strictly semistable. This follows from \cite[Prop 4.6]{NR}, as given $F_1K^*\hookrightarrow F_0\twoheadrightarrow S$ with $F_1$ (and/or $S$) strictly semistable, there exists an irreducible family in $\N_X(n_1,d_1)\times \N_X(n_0-n_1,d_0-d_1+2n_1(g-1))$ with special point $(F_1,S)$ and stable generic point. Since $\delta>n_1(n_0-n_1)(g-1)$, it follows that $\mu(F_1K^*)<\mu(S)$, and so $H^0(S^*F_1K^*)=0$. The latter means that $(F_1,F_0)$ are in the closure of $\pi(\PP(\mathcal{E}xt^1))$. 

For the dimensional computation, if $\delta>n_0n_1(g-1)$, by Riemann--Roch, ${\Zd}'=\N_X(n_1,d_1)\times\N_X(n_0,d_0)$, whence the statement.

 If $(n_0-n_1)n_1(g-1)<\delta\leq n_0n_1(g-1)$, $n_0>n_1$, then for generic $F_0$, when $d_0n_1-(d_1-2(g-1)n_1)n_0> (n_0-n_1)n_1(g-1)$ (namely, for non maximal Toledo by Lemma \ref{lm:toledo_delta}), there always exists a subbundle of the form $F_1K^*$ \cite[Prop. 1.5]{RT}. By \cite[Theorem 0.3]{RT}, for generic $F_0$, the space of bundles with a section $F_1K^*\longrightarrow F_0$ is a $\delta-n_1(n_0-n_1)(g-1)$ dimensional variety. 
 Putting all this together one obtains the result.
\end{proof}
\begin{lemma}\label{lm:sheaf_n0>n1} Let  $\delta\geq n_1(n_0-n_1)(g-1)$ with strict inequality if $n_0\neq n_1$. 

There exist rational projective bundles 
\[
\xymatrix{\PP(\Hh_{\delta})\ar@{-->}[r]&\Zd},\qquad \xymatrix{\PP({\mathcal{\Ee}}_\delta)\ar@{-->}[r] &\Zd},\qquad \xymatrix{\PP({\mathcal{H}}_\delta\oplus \Ee_\delta)\ar@{-->}[r] &\Zd}
\]
with fibers 
\[
\PP(\Hh_\delta)_{(F_0,F_1)}=\PP(H^0(F_1^*F_0K)),\qquad \PP(\Ee_\delta)_{(F_0,F_1)}=\PP(H^1(F_1^*F_0)),\]\[
\PP(\Hh_\delta\oplus\Ee_\delta)_{(F_0,F_1)}=\PP(H^0(F_1^*F_0K)\oplus H^1(F_1^*F_0)).
\]
\end{lemma}
\begin{proof}
By irreducibility of $\Zd$ and density  of $Z^s_{\ol{n},\ol{d}}:=\Zd\cap\N_X^s(n_1,d_1)\times \N_X^s(n_0,d_0)$ for $\delta\geq n_1(n_0-n_1)(g-1)$ with strict inequality if $n_0\neq n_1$ (see Lemma \ref{lm:dimZ}), it is enough to construct the bundles over $Z^s_{\ol{n},\ol{d}}$. Now, over $Z^s_{\ol{n},\ol{d}}$, we may consider the universal extension sheaf \cite{Lange}, yielding $\PP(\Ee_\delta)$  as in Lemma \ref{lm:sheaf_n0=n1}.

To construct $\Hh_\delta$, we may consider the restriction of the universal family of extensions on $\N^s_X(n_0,d_0)\times\N^s_X(n_1,d_1)$  to $\Zd$. By Serre duality, dualisation and  projectivisation yield $\PP(\Hh_\delta)$.

Alternatively, one may consider the \'etale local universal bundles $\Ff_i$ on $X\times\N_X(n_i,d_i)$. Over $X\times Z^s$ the bundles $\PP(\Ff_1^*\boxtimes\Ff_0)$ and $\PP(\Ff_1^*\boxtimes\Ff_0\boxtimes K)$ are well defined and we may argue as in Lemma \ref{lm:sheaf_n0=n1}

The existence of $\PP(\Hh_\delta\oplus\Ee_\delta)$ follows in a similar fashion.
\end{proof}
The following is well known, we sketch a proof for the reader's convenience.
\begin{lemma}\label{lm:irred_comp_from_fixed}
   If the component ${\Cd}$ contains a point with underlying stable bundle, then its underlying reduced scheme ${\Cd}^{red}$ is birational to an affine bundle over $\Fix^{red}$ with fiber $H^1(F_1^*F_0)$ over a stable $(F_0\oplus F_1,\varphi)\in\Fix$.
\end{lemma}
\begin{proof}
     By \cite[Prop. 3.11]{HH}, $\lim_{t\to\infty} t\cdot(E,\phi)=(F_0\oplus F_1,\varphi)$ if and only if $E$ underlies an extension $e\in H^1(F_1^*F_0)$ and $\phi(F_0)=0$, $\Im(\phi)\subset F_0K$ and $\phi$ induces $\varphi$ on the graded object. We may assume $E$ is stable, as this condition is non empty by assumption, hence generic. Likewise, the generic extension in $H^1(F_1^*F_0)$ has underlying stable bundle. Thus, identifying $e$ with a  pair $(i,\pi)$ where $i:F_0\hookrightarrow E$ is the inclusion and $\pi:E\twoheadrightarrow F_1$ the projection, the map
     \[
     e\mapsto (E,i\circ\varphi\circ\pi)
     \]
     defines a birational map from $H^1(F_1^*F_0)$ to the upward flow to $(F_0\oplus F_1,\varphi)$ with inverse given as follows: if $(E,\phi)$ flows up to $(F_0\oplus F_1,\varphi)$, then $(E,\phi)$, from the iterated kernel filtration, it can be assigned an extension class $e'\in H^1((F_1')^*F_0')$ where $F_0':=\Ker\phi\hookrightarrow E\twoheadrightarrow E/\Ker\phi=:F_1'$, which, by \cite[Prop. 3.11]{HH}, satisfies $F_0'\subset F_0$ and $F_1\subset F_1'$. Thus, $F_1'^*F_0'\hookrightarrow F_1^*F_0$, and the induced map $H^1(F_1'^*F_0')\longrightarrow H^1(F_1^*F_0)$ concludes the construction of the inverse. 
     
     We note the map is $\CC^\times$ equivariant. Indeed, if under the aforementioned map we have $e\mapsto (E,i\circ\varphi\circ \pi)$, then $t\cdot e$ is identified with the pair $(t\cdot i,\pi)$ by the extension isomorphism
     \begin{equation}
         \label{eq:ext_iso}
         \xymatrix{F_0\ar[r]^i\ar@{=}[d]&E\ar[r]^\pi\ar[d]^{t\cdot Id}&F_1\ar@{=}[d]\\
F_0\ar[r]_{t\cdot i}&E\ar[r]_{\pi}&F_1,}
     \end{equation}
     so that $t\cdot e\mapsto (E,t\cdot i\circ\varphi\circ \pi)=t\cdot(E, i\circ\varphi\circ \pi)$. Thus, since no other automorphisms of the Higgs bundle exist, the fiber is the whole affine space. Since the upward flow is an affine fibration, this proves that the fibers are as stated over the dense open set of fixed points such that there exists a stable extension in $H^1(F_1^*F_0)$.
     \end{proof}
\begin{theorem}\label{thm:description_fixed_points}
Let $n_0\geq n_1$ with $n_0+n_1=n$, and let $d_0+d_1=d$. Let $\delta$ be as in \eqref{eq:def_delta}. 
 Then $\Fix\subset\M_X(n,d)$ is labelled by $\delta$, and is non-empty if and only if $\delta\in\IC$.

Moreover:
\begin{enumerate}
    \item\label{it:comp_nc} The corresponding reduced component of the nilpotent cone is birational to the  rational bundle 
\[
\xymatrix{\PP(\Ee_\delta\oplus \Hh_\delta)\ar@{-->}[r]&\Zd}
\]
where all sheaves are defined in Lemma \ref{lm:sheaf_n0>n1} and $\Zd$ is defined in Lemma \ref{lm:dimZ}. 

\item\label{it:comp_fp} The scheme  $\Fix^{red}$ is birational to $\PP(\Hh_\delta)$. In particular
\begin{equation}
    \label{eq:dim_Fix}
\dim\Fix=
     (n_0^2+n_1^2-n_0n_1)(g-1)+\delta+1. 
\end{equation}

 \item\label{it:comp_fp_n0>n1} If $n_0>n_1$, the component $\Fix^{red}$ is birational to the rational projective bundle 
\[
\xymatrix{\PP(\tilde{\Ee}_{\delta})\ar@{-->}[r]&\tilde{\Zd}}\]
where
\[
\tilde{\Zd}= \N_X(n_1,d_1)\times \N_X(n_0-n_1,d_0-d_1+2n_1(g-1)) 
\]
and $\PP(\tilde{\Ee_\delta})\longrightarrow \tilde{\Zd}$ has fiber 
\[
\PP(\tilde{\Ee}_\delta)|_{(F_1,S)}=
     \PP(H^1(S^*F_1K^*)).
 \]

\item\label{it:comp_fp_n0=n1} If $n_0=n_1$, then ${\Fd}^{red}$ is birational to the rational bundle
\[\xymatrix{\PP(\tilde{\Hh}_{\delta})/\PP(\Aut(\Oo_\Sym))\ar@{-->}[r]&\Sym^{\frac{\delta}{n_0}}X\times\N_X(n_0,d_0)}\]
where $\tilde{\Hh}_{\delta}$ is defined  in Lemma \ref{lm:sheaf_n0=n1}. 
\end{enumerate}
\end{theorem}
Before tackling the proof of Theorem \ref{thm:description_fixed_points}, we state an immediate corollary, which follows from Theorem \ref{thm:description_fixed_points} via the isomorphism $\M_X(n,d)\longrightarrow \M_X(n,-d)$ given by dualisation on the bundle and transpose on the Higgs field. 
\begin{corollary}\label{cor:comps_n_0<n_1}
    Let $n_0<n_1$.   
    Consider the range
\begin{equation}\label{eq:range_delta*}
       n_0(n_1-n_0)(g-1)<\delta< 2n_0n_1(g-1) \qquad 
\delta\equiv -n_1d+2n_1n_0(g-1)\mod n.
\end{equation}
Then the fixed point components $\Fix\subset\M_X(n,d)$ are labelled by $\delta$, and are non-empty if and only if $\delta$ satisfies \eqref{eq:range_delta*}.

Moreover:
\begin{enumerate}
    \item The corresponding reduced component of the nilpotent cone is birational to the  rational bundle 
\[
\xymatrix{\PP(\Ee_{\delta}\oplus \Hh_{\delta})\ar@{-->}[r]&\Zd}
\]
where $\Zd$ is defined in Lemma \ref{lm:dimZ} and all sheaves are defined in Lemma \ref{lm:sheaf_n0>n1}.

\item The scheme  $\Fix^{red}$ is birational to $\PP(\Hh_{\delta})$. In particular
\[\dim\Fix=     (n_0^2+n_1^2-n_0n_1)(g-1)+\delta+1 
\]

 \item The component $\Fix^{red}$ is birational to the rational projective bundle 
\[
\xymatrix{\PP(\tilde{\Ee}_{\delta}^*)\ar@{-->}[r]&\tilde{\Zd}^*}\]
where
\[
\tilde{\Zd}^*= \N_X(n_0,d_0)\times \N_X(n_1-n_0,d_1-d_0-2n_0(g-1)) 
\]
and $\PP(\tilde{\Ee_{\delta}^*})\longrightarrow \tilde{\Zd}$ has fiber 
\[
\PP(\tilde{\Ee}_{\delta}^*)|_{(F_0,S)}=
     \PP(H^1((F_0K)^*S)).
 \]
\end{enumerate}
\end{corollary}
\begin{proof}
    Let $n_i^*=n_{j}$, ${d}_i^*=-d_{j}$ where $\{i,j\}=\{0,1\}$. Define $\delta^*=d_0^*n_1^*-d_1^*n_0^*+2n_0^*n_1^*(g-1)$.  Then all the statements follow by setting $\mathbf{F}{\ol{n}, \ol{d}}\ni(F_0\oplus F_1,\varphi)=(F_1^*\oplus F_0^*,{}^t\varphi)^*$ and applying Theorem \ref{thm:description_fixed_points} with $(F_1^*\oplus F_0^*,{}^t\varphi)\in\mathbf{F}_{\ol{n}^*, \ol{d}^*}$.
\end{proof}
\begin{remark}
   We note that 
   \[\delta^*=\delta(\ol{n}^*, \ol{d}^*)=d_0^*n_1^*-d_1^*n_0^*+2n_0^*n_1^*(g-1)=\delta(\ol{n}, \ol{d}).\]
   Thus, there is no contradiction between this result and Lemma \ref{lm:toledo_delta}.
\end{remark}
\begin{proof}[Proof of Theorem \ref{thm:description_fixed_points}]
Necessity of non strict bounds follows from Lemma \ref{lm:toledo_delta}.  On the other hand, for $\delta=2n_0n_1(g-1)$, similarly to how we argued in Lemma \ref{lm:toledo_delta}, the  fixed point is of type $n$, and not $(n_0,n_1)$.  Finally, for $n_0>n_1$, and minimal $\delta$,  arguing as in Lemma \ref{lm:toledo_delta}, the fixed point is $S$-equivalent to $(S\oplus F_1K^*\oplus F_1,1_{F_1})$, which is of type $(n_0-n_1,n_1,n_1)$. Indeed, $\mu(F_1K^*)<\mu(S)=\mu(F_1K^*\oplus F_1)=\mu(E)<\mu(F_1)$, so there exists $\beta\in H^0((F_1K^*)^*SK)\neq 0$ and $(S\oplus F_1K^*\oplus F_1,1_{F_1})$ appears as a limit at $0$ of the family $(S\oplus F_1K^*\oplus F_1,t\beta\oplus 1_{F_1})_{t\in\CC}$.

Necessity of the congruence is clear from the expression 
\[
\delta=n_1 d_0-n_0(d-d_0)+2n_0n_1(g-1)=nd_0-n_0d+2n_0n_1(g-1).
\]
We first prove the statement assuming the following:
\begin{assumption}\label{assumption}
    There exists a Higgs bundle $(F_0\oplus F_1,\varphi)\in\Fix$ such that $\varphi$ has maximal rank, and $F_i$ are stable. 
\end{assumption} 
Maximality of the rank means that $\varphi$ is injective and its image is saturated if $n_0>n_1)$ or that $\mathrm{rk} (\varphi_p)=n_1-1$ at each point $p\in\supp(F_0/\Im(\varphi))$ if $n_0=n_1$. Assumption \ref{assumption} will be proven to hold for every component once items \eqref{it:comp_nc}--\eqref{it:comp_fp_n0=n1} have been checked to hold under it.

Sufficiency of the bounds combined with the congruence, and uniqueness of the components for the given invariants follow from the construction below.

Let us prove \eqref{it:comp_nc}. By \cite[Theorem 0.1]{RT}, when  $n_0n_1(g-1)\leq \delta< 2n_0n_1(g-1)$, there exist stable extensions in $H^1(F_1^*F_0)$ for general $F_0, F_1$. Thus, there exists a rational map 
\begin{equation}
    \label{eq:birat_P_C}
    \xymatrix{\PP(\Ee_\delta\oplus\Hh_\delta)\ar@{-->}[r]&{\Cd}\\
(F_0\stackrel{i}{\hookrightarrow} E\stackrel{\pi}{\twoheadrightarrow} F_1,\varphi)
\ar@{|->}[r]& (E,\overbrace{i\circ \varphi \circ\pi)}^{\phi}}
\end{equation}
for some irreducible component ${\Cd}$ of the nilpotent cone with the given invariants. Note that the map is well defined, as the $\CC^\times$ action on the left hand side of \eqref{eq:birat_P_C} acts on the extensions with positive weight and on the Higgs field with negative weight, so $\CC^\times$ orbits of pairs $((i,\pi),\varphi)$ of an extension and a Higgs field are mapped to a unique point. Moreover, the map
 is generically injective. Indeed, if it weren't injective, then
 \[
 (E,i\circ\varphi\circ \pi)\cong (E',i'\circ\varphi\circ \pi') 
 \] 
 Let $f$ be the isomorphism $f\,:\, E\cong E'$ identifying both points.  Since the rank is maximal, the automorphism sends $F_0$ to $F_0$ and $F_1$ to $F_1$, which are stable, and so the induced automorphisms $t_i Id_{F_i}$ are scalar. Now, the extension $(i,t_1\pi)$ is equivalent to $e$, and the extension $(t_0i',\pi')$ is equivalent to $e'$, and $f$ is an isomorphism of extensions $f:(i,t_1\pi)\cong (t_0i',\pi')$. Thus, injectivity follows.
 
 As a consequence, we have 
\[
\dim(\PP(\Ee_\delta\oplus \Hh_\delta)\leq n^2(g-1)+1.
\]
Now, we have
\begin{eqnarray}
    \label{eq:computation_4_BN}
\dim(\PP(\Ee_\delta\oplus \Hh_\delta))\geq -\chi(F_1^*F_0)+\chi(F_1^*F_0K)-1+\dim(\Zd)=
\\\nonumber
\stackrel{\tiny Lm.\ref{lm:dimZ}}{=}3n_0n_1(g-1)-\delta+\delta-n_0n_1(g-1)+(n_0^2+n_1^2)(g-1)+2-1=\\\nonumber=n^2(g-1)+1,   
\end{eqnarray}
where we use $\mu(F_0)<\mu(F_1)$ when $\delta\neq 2n_0n_1(g-1)$.
Note that for any $\Ee=(F_0\oplus F_1,\varphi)\in\Fix$ and any  extension $e\in H^1(F_1^*F_0)$, the underlying bundle $E$ admits a nilpotent Higgs field 
\[
\phi:E\twoheadrightarrow F_1\stackrel{\varphi}{\longrightarrow} F_0K\hookrightarrow EK
\]
that makes it semistable.

A similar argument proves that for $n_1(n_0-n_1)(g-1)\leq \delta\leq n_0n_1(g-1)$ with strict lower bound for $n_0>n_1$, 
\begin{eqnarray}
    \label{eq:computation_4_BN2}
\dim(\PP(\Ee_\delta\oplus \Hh_\delta)\geq -\chi(F_1^*F_0)+1+\dim(\Zd)-1
\\\nonumber
\stackrel{\tiny Lm.\ref{lm:dimZ}}{=}3n_0n_1(g-1)-\delta+1+\delta+1+(g-1)((n_0-n_1)^2+n_0n_1)-1=\\\nonumber=n^2(g-1)+1.
\end{eqnarray}
 To prove \eqref{eq:birat_P_C} is well defined for $n_1(n_0-n_1)(g-1)\leq\delta< n_0n_1(g-1)$ with strict inequality for $n_0>n_1$, we need to take a closer look at $\Zd$. By \cite[Theorem 0.3]{RT}, the generic point $(F_1,F_0)\in \Zd$ satisfies that the general section $\varphi:F_1\longrightarrow F_0K$ is injective and saturated (if $n_0>n_1$) or with maximal rank at the support of $F_0K/F_1$ (if $n_0=n_1$). The condition on $\delta$ implies that $(F_0\oplus F_1,\varphi)$ is stable as a Higgs bundle. Indeed, by the upper bound on $\delta$, $\mu(F_0)\leq \mu(F_1)-(g-1)$,  with equality only if $\delta=n_1(n_0-n_1)(g-1)$. On the other hand, the lower bound of $\delta$ implies that $\mu(F_1\oplus F_1K^*)=\mu(F_1)-g+1\leq\mu(E)$, with equality only if $\delta=n_1(n_0-n_1)(g-1)$. Assuming there exists a Higgs field with maximal rank, any $\varphi$-invariant subbundle is of the form $N_0\oplus N_1$ with $N_i\subset F_i,\, N_1K^*\subset N_0$. Now, by semistability of $F_i$, we have, letting $r_i=\rk(N_i)$,
 \[
 \mu(N_0\oplus N_1)=\frac{\mu(N_0)r_0+\mu(N_1)r_1}{r_0+r_1}\leq \frac{\mu(F_0)r_0+\mu(F_1)r_1}{r_0+r_1}
 \]\[\stackrel{\delta<n_0n_1(g-1)}{<}\mu(F_1)-\frac{n_0n_1(g-1)r_0}{r_0+r_1}
 \leq \mu(F_1)-g+1\leq\mu(E).\]
 Thus, for a general $e\in H^1(F_1^*F_0)$, the image of $(e,\varphi)$ under \eqref{eq:birat_P_C} must also be semistable, and the remaining arguments go through.

Item \eqref{it:comp_fp} is now clear, provided that the rational map be well defined therein. This is guaranteed by the value of $\delta$. Thus,  to prove \eqref{eq:dim_Fix}, we subtract 
\[
\dim\Fix=n^2(g-1)+1-h^1(F_1^*F_0)=(n_0^2+n_1^2-n_0n_1)(g-1)+\delta+1\]
where we have used that $h^0(F_1^*F_0)=0$ generically if $\delta<n_1n_0(g-1)$, and that $\delta>n_1(n_0-n_1)(g-1)$ unless $n_1=n_0$. Thus, \eqref{it:comp_fp} follows.

To see \eqref{it:comp_fp_n0>n1}, recall that $  \delta>n_1(n_0-n_1)(g-1)$. We note that  the existence of $\tilde{\Ee}_\delta$ follows from \cite{Lange}. By irreducibility of all components involved, it is enough to find an injective rational map from one to the other and prove that dimensions match. Now, consider the map
\begin{equation}
    \label{eq:birat_tEe}
\xymatrix{\PP(\tilde{\Ee}_\delta)\ar@{-->}[r]&\Fix^{red}\\
(F_1, S, f_0)\ar@{|->}[r]&\left(F_0\oplus F_1,\varphi:F_1K^*\stackrel{i}{\hookrightarrow}F_0)\right).
}
\end{equation}
In the above $f_0$ denotes an extension class $F_1K^*\stackrel{i}{\hookrightarrow}F_0\stackrel{\pi}{\hookrightarrow}S$ and $F_0$ the underlying vector bundle. 
It is well defined as  by \cite[Theorem 0.3]{RT} the twisted Brill--Noether locus $BN^0_{n_1,d_1-2n_1(g-1)}(F_0)$ (see \eqref{eq:twisted_BN}) is non empty and the general $F_1K^*\in BN^0(F_0)$ is of maximal rank. Moreover, the map \eqref{eq:birat_tEe} is injective, as automorphisms are scalar on both sides by stability of the Higgs bundle for $\delta>n_1(n_0-n_1)(g-1)$. 

To check the equality of dimensions, at a general point 
\[
\dim_{(F_1,S,f_0)} \PP(\tilde{\Ee}_\delta)=
h^1(S^*F_1K^*)+(n_1^2+(n_0-n_1)^2)(g-1)+2-1\]
Since $\mu(S)>\mu(F_1K^*)$ and both bundles are generally stable (so that $H^0(S^*F_1K^*)=0$ for general pairs), by Riemann--Roch
\begin{equation}\label{eq:h1(S*F_1K*)}
h^1(S^*F_1K^*)=-\left(\deg(S^*F_1K^*)-n_1(n_0-n_1)(g-1)\right)=\delta+n_1(n_0-n_1)(g-1),    
\end{equation}
and so
\[
\dim_{(F_1,S,f_0)}\PP(\tilde{\Ee}_\delta)\stackrel{\tiny\eqref{eq:dim_Fix}}{=}\dim\Fix.
\]

For item \eqref{it:comp_fp_n0=n1}, when $n_0=n_1$ $F_1K^*\hookrightarrow F_0$, has torsion cokernel $\mathcal{T}$. There is a map
\begin{equation}
    \label{eq:birat_n0=n1_2}
\xymatrix{\PP(\tilde{\Hh}_\delta)\ar@{-->}[r]&\Fix^{red}\\
(D,F_0,\CC^\times\cdot f)\ar@{|->}[r]&\left(F_0\oplus\ker(f)K, \varphi:\Ker(\CC^\times\cdot f){\hookrightarrow}F_0\right)
}
\end{equation}
which is well defined over points for which $f$ is surjective and hits a component with maximal rank Higgs fields. Note that this map is injective up to the action of a torus $\mathbb{T}:=\Aut(\Oo_D)/\CC^\times \cong (\CC^\times)^{\frac{\delta}{n_0}-1}$, as the kernel of $f\in H^0(F_0^*\Oo_D)$ is well defined up to the action of $\Aut(\Oo_D)$. 

The above arguments prove that there exists a component satifying Assumption \ref{assumption}. We claim that such a component is  uniquely determined by the invariants $(\ol{n},\ol{d})$, or equivalently, $(\ol{n},d_1,\delta)$. Indeed, let $\bfC$ be any component satisfying Assumption \ref{assumption} with fixed invariants, and let ${\Cd}$ be the one admitting a birational map as in \eqref{eq:birat_P_C}. Then, $\bfC$ admits a rational map to $\Zd$. It follows that a subset of the fixed point component of ${\Cd}$ equals the fixed points of $\bfC$. But then these must be dense, as otherwise the dimensional count involved in the proof of \eqref{it:comp_nc}--\eqref{it:comp_fp_n0=n1} would fail.

 We next prove via Brill--Noether theoretic techniques that every component ${\Cd}$ satisfies Assumption \ref{assumption}. The proof is by dimensional count, and based on the constructions within \cite{RT}, so we focus on proving that components have generically maximal rank Higgs field, generic stability of $F_i$ following in a similar fashion.  

In what follows, we assume that $F_i$ are generically stable and prove that non maximality of the rank of the Higgs field yields lower dimensional components. The proof implies that further assuming non stability of the fixed points yields to yet a lower dimensional scheme. 

Now, there are three (non exclusive) options for the  rank not to be maximal: either $\rk(\Ker(\phi))>n_0$, or $F_1$ is not saturated in $F_0$ (when $n_0>n_1$), or $D=\supp(F_0/F_1K^*)$ is never reduced (when $n_0=n_1$). 

\underline{Case 1:} $n_0=n_1$ and $D=\supp(F_0/F_1K^*)$ is never a reduced divisor. Then, there is a rational map
\begin{equation}
    \label{eq:birat_n0=n1}
\xymatrix{\Fix^{red}\ar@{-->}[r]&\PP(\tilde{\Hh}_\delta)\\
(F_0\oplus F_1, \varphi)\ar@{|->}[r]&\left(\mathrm{supp}\left(F_0/F_1K^*\right),F_0,f:F_0{\twoheadrightarrow}F_0/F_1K^*\right)
}
\end{equation}
where $\PP(\tilde{\Hh}_\delta)$ is as in Lemma \ref{lm:sheaf_n0=n1} restricted to the larger diagonal of $\Sym^{\delta/n_0}(X)$. Note that  $\Aut(\mathcal{T})$ contains as a strict subgroup $(\CC^\times)^{\delta/n_0}$. Indeed, since $\deg(\mathcal{T})=\delta/n_0$, as an $\Oo_D$ module, $\mathcal{T}:=F_0/F_1K^*$ is still of rank one. Equivalently, letting $D^r$ denote the reduced divisor underlying $D$, the ranks at the points of $D^r$ equal the multiplicities in $D$. Then $h^0(F_0^*\otimes_{\Oo_D} \mathcal{T})=h^0(F_0^*\otimes_{\Oo_{D^r}} \mathcal{T})=\delta$. Thus
\[
\dim\Fix\leq \PP(\tilde{\Hh}_\delta)/\PP(\Aut(\mathcal{T}))< n_0^2(g-1)+1+s+h^0(F_0^*\otimes \mathcal{T})-1-(\frac{\delta}{n_0}-1)=\]\[
=n_0^2(g-1)+\delta-(\frac{\delta}{n_0}-1-s),
\]
where $s$ is the generic degree of the reduced divisor $D^r$.  The above is strictly smaller than \eqref{eq:dim_Fix}. Since the dimension $h^1(F_1^*F_0)$ does not depend on the rank of the Higgs field, it follows from Lemma \ref{lm:irred_comp_from_fixed} that 
\[
\dim{\Cd}< n_0^2(g-1)+\delta-(\frac{\delta}{n_0}-1-s)+3n_0^2-\delta-1<4n_0^2(g-1)+1. 
\]

\underline{Case 2:} $n_0>n_1$ and generically $\rk(\Ker(\phi))=r>n_0$. Then, for fixed points it must also be $\rk(\Ker(\phi))>n_0$. Then,  there exists  rational map \[\xymatrix{\Fix\ar@{-->}[r]&\PP(\tilde{\Ee}_0)\oplus\PP(\tilde{\Ee}_1)}\] where $\xymatrix{\PP(\tilde{\Ee}_i)\ar@{-->}[r]&\N_X(n_i-r+n_0,d_i-d)\times \N_X(r-n_0,d)}$  has  fiber over $(A,B)$ equal to $\PP(H^1(B^*A))$ if $i=0$ and $\PP(H^1(A^*B))$ if $i=1$. The image of this map hits the restriction to the diagonal $\N_X(n_0-r+n_0,d_0-d)\times \N_X(r-n_0,d)\times \N_X(n_1-r+n_0,d_1-d)\longrightarrow \N_X(n_0-r+n_0,d_0-d) \times \N_X(r-n_0,d) \times \N_X(r-n_0,d)  \times \N_X(n_1-r+n_0,d_1-d)$. A dimensional computation shows that the dimension of these fixed points is strictly less than \eqref{eq:dim_Fix}, and so by Lemma \ref{lm:irred_comp_from_fixed} the total dimension of the component is strictly smaller than half the dimension, a contradiction.

\underline{Case 3:} $n_0>n_1$, $\Im(\varphi) \subset F_0$ never saturated. Then, by Case 2 above we may assume injectivity of $\varphi$. By assumption, all points in this component satisfy that $\varphi\in H^0(\ol{F}_1^*F_0K)$ where $\ol{F}_1$ denotes the saturation of $F_1$ in $F_0K$. Then, let  $\xymatrix{\PP(\tilde{\Ee}_\delta')\ar@{-->}[r]&\N_X(n_1,\ol{d}_1)\times \N_X(n_0-n_1,d_0-\ol{d}_1+2g-2)}$ be as in Lemma \ref{lm:sheaf_n0>n1}, but with $\ol{d}_1$ the generic degree of the saturation of $\Im(\varphi)K^*\subset F_0$. Then, there is a rational map 
\[
\xymatrix{
\Fix\ar@{-->}[r]&\PP(\tilde{\Ee}_\delta')
\\
(F_0\oplus F_1,\varphi)\ar@{|->}[r]& (\ol{F}_1, F_0/\ol{F}_1K^*, f_0)
}
\]
which is well defined and injective as $\varphi\in H^0(\ol{F}_1^*F_0K)$. But $h^1(\ol{F}_1K^*S^*)$ is strictly smaller than \eqref{eq:h1(S*F_1K*)} (by generic stability of $\ol{F}_1$ and $S:=F_0/\ol{F}_1K^*$, so again we have fixed points of dimension strictly smaller than necessary. 
\end{proof}
\begin{remark}
   Uniqueness of the components for a given  $\delta$ also follows from \cite[Corollary 3.3]{Bozec}. This result does not require the hypothesis that underlying stable bundles exist in the component. Our arguments provide a new elementary proof in terms of Brill--Noether theory. It should be noted that, conversely, the results from Higgs bundle theory can be used to recover well known results in Brill--Noether theory. In particular, non emptiness and dimensional computations of twisted Brill--Noether loci follow from the classification of irreducible components of the nilpotent cone. Another example of this interplay is Corollary \ref{cor:BN}. 
\end{remark}
\begin{lemma}\label{lm:dominance}
The rational projections $\xymatrix{\Fix\ar@{-->}[r]&\N_X(n_i,d_i)}$ are dominant  for $\delta\geq n_1(n_0-n_1)(g-1)$, where equality is only considered if $n_0=n_1$.
\end{lemma}
\begin{proof}If $n_0>n_1$, dominance of $\Fix$ over $\N_X(n_0,d_0)$ for non minimal $\delta$ follows from \cite[Prop. 1.11]{RT}. Dominance over  $\N_X(n_1,d_1)$ is a consequence of the description in terms of extensions in Theorem \ref{thm:description_fixed_points}. Similarly, when $n_0=n_1$, dominance of $\Fix$ over $\N_X(n_0,d_0)$ is clear by construction. Dominance over $\N_X(n_1,d_1)$ follows by duality, or from the fact that elementary transformations define a birational morphism between moduli spaces of different degrees. 
\end{proof}

The following well known result in Brill--Noether theory follows from Theorem \ref{thm:description_fixed_points} \eqref{it:comp_nc}.
\begin{corollary}\label{cor:BN}
    Let $n_0\geq n_1$, and let $(F_0,F_1)\in \Zd$ be general. Then
    \[
    h^0(F_1^*F_0K)=\left\{\begin{array}{ll}
  \delta-n_0n_1(g-1)       &  \textrm{ if }\delta\in\IC\setminus\IW  \\
1         & \textrm{ if }\delta\in\IW. 
    \end{array}\right.
    \]
    In particular, let 
    \[BN^0(n_0n_1,\delta):=\{E\in \N_X(n_0n_1,\delta)\,:\, h^0(E)>0\}\] 
    be the Brill--Noether locus
    and recall the twisted Brill--Noether locus $BN^0_{n_1,d_1-2n_1(g-1)}(F_0)$ defined in Equation \eqref{eq:twisted_BN}. Then, for a general $F_0$, tensorisation defines a map \[BN^0_{n_1,d_1-2n_1(g-1)}(F_0)\longrightarrow BN^0(n_0n_1,\delta)\qquad F_1K^*\mapsto F_1^*F_0K\]
    intersecting the open stratum of $BN^0(n_0n_1,\delta)$ consisting of bundles with a minumum number of sections at a non empty subset.
\end{corollary}
\begin{proof}
    For $\delta>n_0n_1(g-1)$, then
    \[
    h^1(F_1^*F_0)+h^0(F_0F_1^*K)-1\stackrel{\tiny Thm. \ref{thm:description_fixed_points}\eqref{it:comp_nc}}{=}\dim\PP(\Ee_\delta\oplus\Hh_\delta)-\dim \Zd\stackrel{\eqref{eq:computation_4_BN}}{=}\]\[=-\chi(F_1^*F_0)+\chi(F_0F_1^*K)-1.
\]
    Since $h^1(F_1^*F_0)=-\chi(F_1^*F_0)$ if $\delta\neq n_1(n_0-n_1)(g-1)$, then it must be $h^0(F_1^*F_0K)-1=\chi(F_1^*F_0K)-1$. This proves the result for $\delta>n_0n_1(g-1)$.

    Similarly,  for  $\delta\leq n_0n_1(g-1)$, then 
    \[
    h^1(F_1^*F_0)+h^0(F_0F_1^*K)-1\stackrel{\tiny Thm. \ref{thm:description_fixed_points}\eqref{it:comp_nc}}{=}\dim\PP(\Ee_\delta\oplus\Hh_\delta)-\dim \Zd\stackrel{\eqref{eq:computation_4_BN2}}{=}-\chi(F_1^*F_0)+1-1.
    \]
Thus $h^0(F_1^*F_0K)=1$, as $h^1(F_1^*F_0)=-\chi(F_1^*F_0)$. For if $\delta= n_1(n_0-n_1)(g-1)$, $h^1(F_1^*F_0)=h^0(F_1^*F_1K^*)\oplus h^0(F_1^*SK^*)$, which is generally zero by Theorem \ref{thm:description_fixed_points}\eqref{it:comp_fp_n0=n1}.
\end{proof}

\section{Wobbly divisors of type $(n_0,n_1)$}\label{sec:wobbly_divisors}
In this section we prove Drinfeld's conjecture for a particular kind of wobbly vector bundles: those of wobbly-type $(n_0,n_1)$ (cf. Definition \ref{def:wobbly-type}).   This generalises the results in \cite{PalPauly, PPrk3} for the rank two and three cases.

Recall the forgetful morphism
\begin{equation}\label{eq:forgetful}
\M_X(n,d)\stackrel{f}{\dasharrow}\N_X(n,d).
\end{equation}
Define \[\Wd:=\ol{\Im(f|_{\tiny \Cd}^\delta})\subset\N_X(n,d).\] 
We first check that this is non empty for all $\delta\in\IC$.
\begin{lemma}
   \label{lm:comps_intersect_N} Let $\delta\in\IC\setminus\IW\cup\{\delta_m\}$, where $\delta_m=\max\IW$.
   Then $\Cd\cap\N_X(n,d)\neq\emptyset$. 
\end{lemma}
\begin{proof}
    For $\delta> n_0n_1(g-1)$, it was shown in the proof of Theorem \ref{thm:description_fixed_points}. For $\delta\in\IW$, define $s_\delta=2n_0n_1(g-1)\geq n_0n_1(g-1)$. Then, every stable bundle $E$ contains a $F_0$ subbundle of degree $d_0$ and rank $n_0$ \cite[Thm. 0.3]{RT}. We need to show that there exist $F_0, E$ as stated with $H^0\left(E/F_0)^*F_0K\right)\neq 0$.
     By Lemma \ref{lm:decreasing_delta}, this will hold for $\delta_m\in \{n_0n_1(g-1),  n_0n_1(g-1)-1, n_0n_1(g-1)-2\}$. 
\end{proof}
\begin{theorem}
    \label{thm:wobbly_divisors}
 Let $\bfW_{n_0,n_1}:=\bigcup_{\delta\in\IC}\Wd$ be the set of wobbly-type $(n_0,n_1)$ wobbly vector bundles. Then
 \begin{enumerate}
     \item\label{it:pure_codim} $\bfW_{n_0,n_1}$ is of pure codimension one with decomposition into irreducible components $\bfW_{n_0,n_1}=\bigcup_{\delta\in I}\Wd$ with $\delta_m\in I\subseteq\IW$ defined by
     \[
     I=\{\delta\in\IW\,:\,\Cd\cap\N_X(n,d)\neq\emptyset.\}.
     \]
     \item\label{it:Wdelta_Wdelta-n}  Given $\delta\in\IC\setminus\IW$, then
$\Wd\subsetneq\bfW^{\delta-n}_{(n_0,n_1)}\subsetneq \bfW^{\delta_m}_{(n_0,n_1)}$, where $\delta_m:=\max\IW$.
 \end{enumerate}
\end{theorem}
\begin{proof}
In what follows, we consider the underlying reduced schemes. The case $n_1=n_0=1$ follows from \cite{PalPauly}, so we may assume $n\geq 3$\footnote{The case $n=3$ is shown in \cite{PPrk3}, but we hereby give a different proof.}. 

The fact that $\delta_m\in I$ follows from Lemma \ref{lm:comps_intersect_N}.

To prove \eqref{it:pure_codim}, we will show that $\Wd$ is of codimension one if and only if it is non-empty and $\delta\in\IW$ and then use \eqref{it:Wdelta_Wdelta-n} to conclude. 

Note that it is enough to consider a dense open set inside $E\in \Wd$. Thus, we may assume $E\in\bfW_\delta^{smooth}\cap \N_X^s(n,d)$ (since the general bundle underlying a point in $\Cd$ is stable by Theorem \ref{thm:description_fixed_points}).

To see the necessity of $\delta\in\IW$, consider the forgetful map $f|${\tiny$_{\bfC^\delta_{\ol{n}}}$}.
This map is the composition
\[
\Cd\hookrightarrow\M_X\stackrel{\delta}{\longrightarrow}h^{-1}(0)^{\CC^\times}\dasharrow \N_X\]
where 
\[
\delta: \M_X\longrightarrow h^{-1}(0)^{\CC^\times}\subset\M_X
\]
is the downward flow $(E,\varphi)\mapsto \lim_{t\to0}(E,t\varphi)$. This is a Zariski locally trivial fibration, with $H^0(\End(E)K)=\delta^{-1}(E)$.  Also, given $\varphi\in H^0(F_1^*F_0K)$ we see that 
\[
df: T_{(E,\varphi)}\M_X\longrightarrow T_E\N_X
\]
is just the composition of 
\[
d\delta: T_{(E,\varphi)}\M\longrightarrow T_{(E,0)}\M_X
\] and the projection onto the first factor of the decomposition $T_{(E,0)}\M_X\cong T_E\N_X\oplus H^0(\End(E)K)$.  

Now, by triviality of the fibration and smoothness of $(E,0)$, then $H^0(\End(E)K)=N_{\N_X}\M_X$ is the normal bundle to $\N_X\subset \M_X$. It follows that $N_{\bfW}\Cd\supset N_{\N_X}\M_X \cap T_{E,0}\Cd=H^0(F_1^*F_0K)$.

Sufficiency of $\delta\in\IW$ follows from Lemma \ref{lm:unique_filtration}, which shows local injectivity of the downward flow map $\delta:\Cd\longrightarrow\Fixd$.
\end{proof}
\begin{lemma}
    \label{lm:unique_filtration}
    Let $\delta\in\IW$. Let $E\in\Wd$ be general. Then, there exists a unique nilpotent $\varphi\in H^0(\End(E)\otimes K)$ such that $(E,\varphi)\in\Cd$.
\end{lemma}
\begin{proof}
    Assume there existed $\varphi\neq\varphi'$. Then, $E$ can be expressed as an extension of bundles with fixed invariants $(n_0,n_1), (d_0,d_1)$ in two ways:
    \[
    F_0=\Ker(\varphi)\hookrightarrow E\twoheadrightarrow F_1'.
    \]  \[
    F_0'=\Ker(\varphi')\hookrightarrow E\twoheadrightarrow F_1'.
    \]
     Since $E$ is generic, we may assume that $(F_0,F_1), \in \Zd$ is also generic. 
     
    Then, there exists an exact diagram
    \[
\xymatrix{&0\ar[d]&0\ar[d]&0\ar[d]&\\
0\ar[r]&A\ar[d]\ar[r]&F_0'\ar[r]\ar[d]&B'\ar[r]\ar[d]&0\\
0\ar[r]&F_0\ar[d]\ar[r]&E\ar[r]\ar[d]&F_1\ar[r]\ar[d]&0\\
0\ar[r]&B\ar[r]\ar[d]&F_1'\ar[d]\ar[r]&M\ar[d]\ar[r]&0\\
&0&0&0&.}
    \]
    Note that by Lemma \ref{lm:dominance}, we may assume that $F_0$ is general inside $\N_X(n_0,d_0)$. Then, $A, B$ can be supposed to be stable by \cite[Proposition 2.6]{NR}. 
    Then
\[\mu(B)-g+1=\mu(B'K^*)<\mu(F_1K^*)<\mu(F_0)<\mu(B')=\mu(B)<\mu(F_1)\stackrel{\delta\leq n_0n_1(g-1)}{\leq}\mu(F_0)+g-1.
\]
This implies that generically  $h^0(F_0^*B)=0$, contradicting the above diagram unless $B=F_0$, in which case it must be $n_0=n_1$ and $M$ torsion. Then
\[
e\in\Ker(H^1(F_1^*F_0)\twoheadrightarrow H^1((F_1')^*F_0))
\]
varies in a subspace of codimension $-\chi((F_1')^*F_0)=-\delta+3n_0^2(g-1)$.

Now, fixing $(F_0,F_1)\in\Zd$,  $F_1'$ varies in a space of dimension $\dim\PP(\Ext^1(M,F_0))=\delta-1$. On the other hand, $M$ varies in a space of dimension at most $\deg(M)=\frac{\delta}{n_0}$. Thus, altogether we have that $E$ varies in a space of codimension at least 
\[
\codim\geq 3n_0^2(g-1)-2\delta-\frac{\delta}{n_0}+1\stackrel{\delta<2n^2_0(g-1)}{\geq}1.
\]
\end{proof}
\begin{lemma}\label{lm:decreasing_delta}
    Let $n\geq 3$, and let $F_0, F_1$ be stable vector bundles of ranks $n_0\geq n_1$. Let $\delta=\deg(F_1^*F_0K)$. Assume that $\delta>n_0n_1(g-1)+1$. Then, there exist vector bundles $F_0'$, $F_1'$ fitting in exact sequences
    \begin{eqnarray}
        \label{eq:Heckes}
    F_0'\hookrightarrow F_0\twoheadrightarrow\Oo_p
    \\
    \nonumber
    F_1\hookrightarrow F_1'\twoheadrightarrow\Oo_q
      \end{eqnarray}
    and a section $\varphi\in H^0(F_1^*F_0K)$ factoring as
    \[
\varphi\,:\,    F_1\hookrightarrow F_1'\stackrel{\varphi'}{\longrightarrow} F_0'K\hookrightarrow F_0K.
    \]
    Moreover, given $F_i$, $F_i'$, $p\in X$ satisfying \eqref{eq:Heckes}, then the set of $E'\in \PP(H^1((F_1')^*F_0'))$ may be identified with $\PP(\Ext^1(\Oo_p,E'_p)\setminus\PP((F_0')_p)\subset\PP^{n-1}$, while the set of $E'\in \PP(H^1(F_1^*F_0)$ defines the closed subset $\PP((F_0')_p)\subset\PP^{n-1}$.
\end{lemma}
\begin{proof}
Under the given hypotheses, we have that 
\[[\cdot]\colon H^0(F_1^*F_0K)\longrightarrow 
H^0(F_1^*F_0K)/H^0(\End(F_1))\] is not equivalently $0$ for dimensional reasons. Then, we may choose $\varphi\in H^0(F_1^*F_0K)$ with $[\varphi]\neq 0$. By Lemma \ref{lm:perp2}, it follows that if $n_0>n_1$, $\rk(\varphi_p)\leq n_1-1<n_0-1$, so we may apply Lemma \ref{lm:factor} to find the desired $F_i'$ and $\varphi'$. If $n_0=n_1$, consider the exact sequence
\[
N\hookrightarrow F_1\stackrel{\varphi}{\longrightarrow} F_0K\twoheadrightarrow SK.
\]
If $\rk(N)>1$, or $N\neq 0$ and $S$ has torsion (e.g., when $\rk(N)=1$ for degree reasons), then we may apply Lemma \ref{lm:factor} to conclude. If $N=0$, $S$ is a torsion sheaf, and we note that $\deg(S)=\delta/n_0>n_0$, so either there exists $p\in\supp(S)$ such that $\dim(S_p)=2$ (in which case we apply Lemma \ref{lm:factor}), or there exist $p\neq q$, $p,q\in\supp(S)$. In the latter case, we may consider Hecke transforms at possibly different points $p, q$, as in Equation \eqref{eq:Heckes}, in such a way that the Higgs field factors as desired.

  Now, consider a stable bundle $E\in \PP(H^1(F_1^*F_0))$. We want to produce a family of extensions inside $\PP(H^1((F_1')^*F_0'))$ degenerating to $E$. Since $E$ is stable by assumption, so will be the generic $E'$ inside the given family. 

  Observe that there are surjections
\[H^1((F_1')^*F_0')\twoheadrightarrow H^1((F_1')^*F_0)\twoheadrightarrow H^1(F_1^*F_0)\]
  We consider $e'\in H^1((F_1')^*F_0')$, $\epsilon\in H^1((F_1')^*F_0)$ with
\[
e'\mapsto\epsilon\mapsto e.
\]
This is always possible as all morphisms involved are surjective. Then, the underlying bundles fit in a long exact sequence
 \[
    \xymatrix{
F_0'\ar@{^(->}[d]\ar@{^(->}[r]&E'\ar@{->>}[r]\ar[d]&F_1'\ar@{=}[d]\\
F_0\ar@{^(->}[r]&\Ee\ar@{->>}[r]&F_1'\\
F_0\ar@{^(->}[r]\ar@{=}[u]&E\ar@{^(->}[u]\ar@{->>}[r]&F_1\ar@{^(->}[u],}
    \]
and $E'\in\bfW_{\delta-n}^{n_0,n_1}$.
Then, since $F_i$ are stable, then $F_i'$ are at worst semistable. Also, since 
\[\frac{\delta}{n_0}=\mu(F_1)-\mu(F_0)>n_0\geq 2,\]
then $\mu(F_0')<\mu(F_1')$, $\mu(F_0)<\mu(F_1)$, and we have a map
    \[
\PP(H^1((F_1')^*)F_0')\hookrightarrow
   \PP(H^1((F_1')^*)F_0)\]
   mapping $E'$ to $\Ee\in \PP(\Ext^1(\Oo_q,E'))\cong\PP^{n-1}$. We want to characterise the bundles $\Ee$ which contain $E\in H^1(F_1^*F_0)$. Let $l_p\in\PP((F_0')_p)$ be such that $l_p\subseteq\Ker( (F_0')_p\longrightarrow (F_0)_p)$. Then, the extensions $\Ee$ such that $F_0\subset \Ee$ are precisely those containing $l_p$. Hence, they are identified with $\PP((F_0)_p')\subset \PP(E_p')\cong\PP(\Ext^1(\Oo_p,E')).$
\end{proof}
\begin{lemma}\label{lm:factor}
    Let $\Ee\in \Fix$ be generic. Assume that there exists $p\in X$ such that
        $\rk(\varphi_p)<n_0-1$. 
   Then, 
    there exist Hecke transforms
     \[
   \tilde{F}_0\stackrel{i_0}{\hookrightarrow }{F}_0\twoheadrightarrow \Oo_p
    \]
    \[
    F_1\stackrel{i_1}{\hookrightarrow}\tilde{F}_1\twoheadrightarrow \Oo_p,
    \]
    and $\tvarphi\in H^0(\tF_1^*\tF_0K)$ such that $\tvarphi=i_0\circ\varphi\circ i_1$.
\end{lemma}
\begin{proof}
  By hypothesis, there exists a plane $\Pi_p\subset (F_0K)_p$ such that
  $\Im(\varphi_p)\subset\Ker((F_0K)_p\longrightarrow \Pi_p)$. 

  Consider $l_p\subset \Pi_p$ a line therein. Then $\Im(\varphi_p)\subset\Ker((F_0K)_p\longrightarrow l_p)$. Thus, we may define $\tilde{F}_0=\Ker(F_0\longrightarrow l_p)$, where abusing notation we identify $(F_0K)_p\cong (F_0)_p$. Moreover, 
  \[
\dim\Ker((F_1)_p\stackrel{\varphi_p}{\longrightarrow}(\tF_0K)_p)\geq 1.\]
Let $m_p\subseteq \Ker((F_1)_p\stackrel{\varphi_p}{\longrightarrow}(\tF_0K)_p)$ be a line.  We may thus define
  \[
\tF_1=(\Ker(F_1^*\longrightarrow m_p^*))^*.
  \]
\end{proof}
\begin{remark}
Let $n_0\geq \frac{3}{2}n_1$. Then $n_0n_1(g-1)+1<3n_1(n_0-n_1)(g-1)$, and so for all $\delta\notin \IW$, $\Wd\subset \bfW_{(n_1,n_1,n_0-n_1)}\cap \bigcap_{\delta'>\delta}\bfW_{(n_0,n_1)}^{\delta'}$ by Theorems \ref{thm:wobbly_comps} and \ref{thm:wobbly_divisors}.
\end{remark}
\begin{lemma}\label{lm:perp2}
    Let $\Ee\in \Fix$ be generic, $\delta>0$. Assume that 
      there exists $p\in X$ such that
        $\rk(\varphi_p)<n_1$. 
    Then, $\varphi\in H^1(\End(F_1))^\perp\cap H^0(F_1^*F_0K)\subset H^0(\End(F_0\oplus F_1)K)$. Equivalently, consider the natural quotient map
    $[\cdot ]\,:\, H^0(F_1^*F_0K)\twoheadrightarrow H^0(F_1^*F_0K)/H^0(\End(F_1)K)$. Then, $[\varphi]\neq 0$.

 Conversely, if $[\varphi]\neq 0$ then $\rk(\varphi_p)<n_1-1$ for some $p\in X$.
\end{lemma}
\begin{proof}
    Assume that $\rk(\varphi_p)<n_1$. Consider the exact sequence
    \[
    N\hookrightarrow F_1\stackrel{\varphi}{\longrightarrow} F_0K\twoheadrightarrow SK,
    \]
    where $N$ may be zero or not.
    Then, there is an exact diagram
    \[
    \xymatrix{
 S^*NK^*\ar@{^(->}[r]\ar@{^(->}[d]&S^*F_1K^*\ar@{^(->}[d]\ar[r]&S^*F_0\ar@{^(->}[d]\ar@{->>}[r]&\End(S)\ar@{^(->}[d]\\
    F_0^*NK^*\ar@{^(->}[r]\ar[d]&F_0^*F_1K^*\ar[d]\ar[r]&\End(F_0)\ar[d]\ar@{->>}[r]&F_0^*S\ar[d]\\
   F_1^*N\ar@{^(->}[r]\ar[d]&\End(F_1)\ar[d]\ar[r]&F_1^*F_0K\ar[d]\ar@{->>}[r]&F_1^*SK\ar[d]\\
   \End(N)\ar@{^(->}[r]&N^*F_1\ar[r]&N^*F_0K\ar@{->>}[r]&N^*SK.
    }
    \]
    By the genericity hypothesis on $\Ee$ we may assume that $F_i, S$ be stable. 
    
    Considering the associated long exact diagram in cohomology
    \[
     \xymatrix{
H^0(S^*NK^*)\ar@{^(->}[r]\ar@{^(->}[d]&\overbrace{H^0(S^*F_1K^*)}^{=0}\ar@{^(->}[d]\ar[r]&\overbrace{H^0(S^*F_0)}^{=0}\ar@{^(->}[d]\ar[r]&\overbrace{H^0(\End(S))}^{=\CC}\ar@{^(->}[d]...\\
  H^0(F_0^*NK^*)\ar@{^(->}[r]\ar[d]&\overbrace{H^0(F_0^*F_1K^*)}^0\ar[d]\ar[r]&H^0(\End(F_0))\ar[d]\ar[r]&H^0(F_0^*S)\ar[d]...\\
H^0(F_1^*N)\ar@{^(->}[r]\ar[d]&H^0(\End(F_1))\ar[d]\ar[r]&H^0(F_1^*F_0K)\ar[d]\ar[r]&H^0(F_1^*SK)...\ar[d]\\
H^0( \End(N))\ar@{^(->}[r]&H^0(N^*F_1)\ar[r]&H^0(N^*F_0K)\ar[r]&H^0(N^*SK)...
    }
    \]
    with suitable modifications if $N=0$.
    Note that $H^0(\End(F_1))=H^1(\End(F_1))^*\cap H^0(F_1^*F_0K)$. Thus, $\varphi\in H^1(\End(F_1))^\perp\cap H^0(F_1^*F_0K)$ if and only if $\varphi$ maps to a non zero element of $H^0(F_1^*SK)$. If $\rk(\varphi_p)<n_1$, it follows that $\varphi=0$ or $[\varphi]\neq 0$, as otherwise $\varphi=\lambda\mathrm{Id}_{F_1}$ and the rank is constantly $n_1$. 
    
    Conversely, if $[\varphi]\neq0$, then if $N\neq 0$ we are done. If $N=0$, then $F_1$ is a subsheaf of $F_0K$ with cokernel supported at those points over which $[\varphi]_p\neq 0$.
\end{proof}

\section{Criteria for wobbliness}\label{sec:criteria}
In this section we give some criteria to decide when a fixed point is wobbly. We start by defining types of wobbly Higgs bundles, dependant of the upward flow from them. 
\begin{definition}\label{def:wobbly-type}
Let $\ol{m}$ be an ordered partition of $n$. A wobbly fixed point $\Ee$ is called of wobbly-type $\ol{m}$ if $\Ee=\lim_{t\to 0}t\cdot(E,\psi)$ and  $\Ee\neq\Ee'=\lim_{t\to \infty}t\cdot(E,\psi)$ is of type $(\ol{m},\ol{d})$ (cf. Definition \ref{def:type_fixed_point}). 
\end{definition}
We will denote by $\HWFix$ the set of points of wobbly-type $\ol{m}$ inside $\Fix$, that is
\begin{equation}
    \label{eq:HWFix}
\HWFix=\{\Ee\in\Fix\colon\exists\Ee'\in\mathbf{F}_{\ol{m}}\cap \ol{\Ee^+}.\}
\end{equation}
It follows that 
\[
\HWFix=\bigcup_{\bfC_{\ol{m}}^{\ol{d}}\cap\Fix\neq0}\HWFixd
\]
where 
\begin{equation}
    \label{eq:HWFixd}
\HWFixd=\{\Ee\in\Fix\colon\exists\Ee'\in\mathbf{F}_{\ol{m}}^{\ol{d}}\cap \ol{\Ee^+}.\}
\end{equation}
\begin{remark}
    The wobbly-type is not uniquely determined. For type $(n)$ wobbly points, any point in the intersection of two irreducible components has at least two wobbly types (cf. \cite{PalPauly} for the rank two case). In the current context, see Remark \ref{rk:wobbly_type}.
\end{remark}
Let us recall a result from \cite{HH} that will appear many times in this section.
\begin{lemma}{\cite[Props. 3.4 \& 3.11]{HH}}\label{prop:3.4_HH}
A Higgs bundle $(E,\psi)$ satisfies $\lim_{t\to 0}t(E,\psi)=(F_0\oplus F_1,\varphi)$ where $(F_0\oplus F_1,\varphi)\in\M_X^s$ if and only if $E$ is an extension of $F_0$ by $F_1$ and $\psi$ induces $\varphi$ under the natural maps
\[
\varphi\,:\,F_1\hookrightarrow{E}\stackrel{\psi}{\longrightarrow} EK\twoheadrightarrow {F}_0K.
\]

Similarly, $\lim_{t\to \infty}t(E,\psi)=(F_0\oplus F_1,\varphi)$ where $(F_0\oplus F_1,\varphi)\in\M_X^s$ if and only if $E$ is an extension of $F_1$ by $F_0$, and $\psi|_{F_0}\equiv 0$ with the induced factorisation of $\psi$ matching $\varphi$, namely:
\[
\psi\,:\,E\twoheadrightarrow F_1\stackrel{{\varphi}}{\longrightarrow} {F}_0K\hookrightarrow EK.
\]
\end{lemma}
Let us start by a preliminary lemma.
\begin{lemma}\label{lm:perp}
Let $F_1,\, F_0$ be stable bundles with $\mu(F_1)\geq\mu(F_0)$, $F_0\not\cong F_1$, and let $e\in H^1(F_0^*F_1)$ be an extension class with underlying vector bundle $E$. Let
\[
Res:H^0(\End(E)\otimes K)\longrightarrow H^0(F_0F_1^*K)
\]
be the natural map. Then, given $\varphi\in H^0(F_0F_1^*K)$, the following are equivalent:
\begin{enumerate}
    \item\label{it:phi_lifts}There is an equality $\varphi=Res(\varphi)$ for some $\varphi\in H^0(\End(E)\otimes K)$.
    \item\label{it:phi_perp_e} The pairing given by Serre duality satisifies $\langle\varphi,e\rangle=0$.
    \item\label{it:phi_lifts_01} The field $\varphi$ lifts to $H^0(F_1^*EK)$ and $H^0(F_0E^*K)$.
    \item\label{it:phi_lifts_0_or_1} The field $\varphi$ lifts to $H^0(F_1^*EK)$ or $H^0(F_0E^*K)$.
\end{enumerate} 
\end{lemma}
\begin{proof}
Consider the exact diagram
\[
\xymatrix{
F_1F_0^*K\ar@{^(->}[r]\ar@{^(->}[d]&F_1{E}^*K\ar@{->>}[r]\ar@{^(->}[d]&\End(F_1)K\ar@{^(->}[d]\\
EF_0^*K\ar@{^(->}[r]\ar@{->>}[d]&\End(E)K\ar@{->>}[r]\ar@{->>}[d]&EF_1^*K\ar@{->>}[d]\\
\End(F_0)K\ar@{^(->}[r]&F_0{E}^*K
\ar@{->>}[r]&F_0F_1^*K.
}
\]
The associated long-exact diagram in cohomology reads
\begin{equation}
    \label{eq:LES_lift}
    \xymatrix{
H^0(F_1F_0^*K)\ar@{^(->}[r]\ar@{^(->}[d]&H^0(F_1{E}^*K)\ar[r]\ar@{^(->}[d]&H^0(\End(F_1)K)\ar[r]^{-\cup e}\ar@{^(->}[d]&H^1(F_1F_0^*K)\ar[d]\ar[r]&\dots\\
H^0(EF_0^*K)\ar@{^(->}[r]\ar[d]&H^0(\End(E)K)\ar[r]\ar[d]&H^0(EF_1^*K)\ar[r]^{-\cup e}\ar[d]&H^1(EF_0^*K)\ar[d]\ar[r]&\dots\\
H^0(\End(F_0)K)\ar@{^(->}[r]\ar[d]_{\cup e}&H^0(F_0{E}^*K)\ar[d]_{\cup e}\ar[r]&H^0(F_0F_1^*K)\ar[r]^{-\cup e}\ar[d]_{\cup e}&H^1(\End(F_0)K)\ar[r]\ar[d]&\dots\\
H^1(F_1F_0^*K)\ar[r]&H^1(F_1{E}^*K)\ar[r]&H^1(\End(F_1)K)\ar[r]&0&\\
}
\end{equation}
Let $\varphi\in H^0(F_0F_1^*K)$. 

It is clear that \eqref{it:phi_lifts} $\Rightarrow$ \eqref{it:phi_lifts_01}$\Rightarrow$ \eqref{it:phi_lifts_0_or_1}.

Let us prove that \eqref{it:phi_lifts_0_or_1} $\Leftrightarrow$ \eqref{it:phi_perp_e} $\Rightarrow$ \eqref{it:phi_lifts}, which finishes the proof. 

Let $\varphi$  lift to $\varphi_1\in H^0(EF_1^*K)$. Note that we have  
\[\langle e, \cdot\rangle\,:\,H^0(F_0F_1^*K)\stackrel{\cup e}{\longrightarrow} H^1(\End(F_1)K)\cong \CC.\] 
Thus, sections lifting to $H^0(EF_1^*K)$ are contained in the hyperplane $e^\perp\subset H^0(F_0F_1^*K)$, where $e^\perp=\{\psi\in H^0(F_0F_1^*K)\,:\,\langle e,\psi\rangle=0\}$; so \eqref{it:phi_perp_e} holds. The same argument proves that if $\varphi$  lifts to $\varphi_2\in H^0(F_0{E}^*K)$, then it must be $\varphi\perp e$ (i.e., $\langle e,\varphi\rangle=0$). Now, assume $e\perp\varphi$. Note that $H^0(EF_1^*K)$ surjects onto ${e}^\perp\subset H^0(F_1^*F_0K)$, as the cokernel of $H^0(EF_1^*K)\longrightarrow H^0(F_0F_1^*K)$ is either $0$ or $H^1(\End(F_1)K)\cong\CC$ (by stability of $F_1$). Hence $e\perp\varphi$ implies \eqref{it:phi_lifts_0_or_1} (and in fact it implies \eqref{it:phi_lifts_01}). Now, continuing the argument  $\varphi_1\cup e\in H^1(EF_0K)$ maps to zero under $H^1(EF_0^*K)\longrightarrow H^1(\End(F_0)K)$. Now, the kernel of this map is $H^1(F_1F_0^*K)\cong H^0(F_1^*F_0)^*=0$ by stability of $F_i$ and the slope inequality. It thus follows that $\varphi_1\cup e=0$, so that $\varphi_1$ lifts to $\varphi\in H^0(\End(E)K)$.
\end{proof}
\begin{proposition}\label{prop:wobbly_n0n1}
     Let $n_1\leq n_0$. A generic fixed point $\Ee=(F_0\oplus F_1,\varphi)\in\Fix$ is of wobbly type $(n_0,n_1)$ if and only if there exists a non zero torsion sheaf $\mathcal{T}$, and elementary transformations
     \begin{eqnarray}
         \label{eq:tF0olF1}
\tilde{F}_0\stackrel{i_0}{\hookrightarrow} F_0\stackrel{f_0}{\twoheadrightarrow} \Tt
     \\\nonumber
     {F}_1\stackrel{i_1}{\hookrightarrow} \ol{F}_1\stackrel{f_1}{\twoheadrightarrow} \Tt
     \end{eqnarray}
     such that $\varphi$ factors as 
\[F_1\stackrel{i_1}{\hookrightarrow}\ol{F_1}\stackrel{\ol{\varphi}}{\longrightarrow }\tilde{F}_0K\stackrel{i_0}{\hookrightarrow} F_0K.\] 
\end{proposition}
\begin{remark}
    To illustrate Proposition \ref{prop:wobbly_n0n1}, let us analyse two simple examples: when $\rk(F_0)=\rk{F_1}=1$, then the criterion is equivalent to $\varphi$ having a zero of order at least two. Indeed, if the latter holds, then we may take
    $\mathcal{T}=\Oo_p$, and  the map $F_0\twoheadrightarrow \Oo_{p}$; it then follows that $\varphi$ factors as stated in Proposition \ref{prop:wobbly_n0n1}. Conversely, it is easy to see that for line bundles the criterion is only realised in terms of vanishing of the Higgs field. In this particular case, Proposition \ref{prop:wobbly_n0n1} follows from  \cite[Theorem 4.16]{HH}. On the other hand, for $\rk(F_0)=2, \rk{F_1}=1$, simple vanishing of $\varphi$ is a generic condition ensuring $(2,1)$-wobbly type. Indeed, there is always an elementary transformation $F_0\longrightarrow \Oo_p$ containing $F_1(p)K^*$ \cite{PPrk3}.
\end{remark}
\begin{proof}
We note that if $(F_0\oplus F_1,\varphi)$ is (semi)stable, then, under the conditions of the statement, so is $(\tilde{F}_0\oplus \ol{F}_1,\ol{\varphi})$. We claim that  $\ol{\varphi}$-invariant subbundles of $\tilde{F}_0\oplus \ol{F}_1$ are of the form $\tilde{N}\oplus \ol{M}$ where $\tilde{N}=N\cap\tilde{F}_0$ for some $N\subset F_0$, $\ol{M}\supset M$ for some $M\subset F_1$, $\varphi(M)\subset NK$ factors as $\varphi:M\subset\ol{M}\stackrel{\ol{\varphi}}{\longrightarrow} \tilde{N}K\subset NK$ and $\deg{M\oplus N}=\deg(\ol{M}\oplus \tilde{N})$. To see this, define $M:=F_1\cap\ol{M}$, so that $\Tt_M:=\ol{M}/M\subset \Tt$ is torsion. Define $N$ to be the saturation of $\tilde{N}$ inside $F_0$ along $\Tt_M$, so that $\tilde{N}:=N\times_{\Tt}\Tt_M$. It is clear that $\varphi$ factors as stated. Also, $\deg ({N})-\deg (\tilde N)=\deg(\ol{M})-\deg(M)$, so that $\mu(\tilde{N}\oplus\ol{M})=\mu({N}\oplus{M})\leq\mu(F_0\oplus F_1)=\mu(\tilde{F}_0\oplus\ol{F}_1)$ by semistability of $F_0\oplus F_1$ and an equality happens if $F_0\oplus F_1$  is strictly semistable. .

Let us first see that the existence of elementary transformations as in \eqref{eq:tF0olF1} satisfying the given condition is sufficient. We claim that there exists an extension $e\in\Ker\left(H^1(\ol{F}_1^*\tilde{F}_0)\stackrel{h^1(i_0)}{\longrightarrow} H^1(\ol{F}_1^*{F}_0)\right)$ whose underlying vector bundle $E$ fits in an exact commutative diagram \begin{equation}
    \label{eq:good_cd}
\xymatrix{
&F_1\ar@{^(->}[d]\ar@{=}[r]&F_1\ar@{^(->}[d]\\
\tilde{F}_0\ar@{=}[d]\ar@{^(->}[r]&{E}\ar@{->>}[d]\ar@{->>}[r]&\ol{F}_1\ar@{->>}[d]^{f_1}\\
\tilde{F}_0\ar@{^(->}[r]&{F}_0\ar@{->>}[r]_{f_0}&\mathcal{T}.
}.
\end{equation}
In that case, we have that the Higgs field \[
\psi\,:\,E\twoheadrightarrow \ol{F}_1\stackrel{\ol{\varphi}}{\longrightarrow} \tilde{F}_0K\hookrightarrow EK
\] 
is nilpotent, and satisfies the conditions of Proposition \ref{prop:3.4_HH}, so
\[\lim_{t\to 0} (E,\psi)=\Ee; \qquad\lim_{t\to \infty} (E,\psi)=\left(\tF_0\oplus\lF_1,\ol{\varphi}\right),\]
and so $\Ee$ is of wobbly-type $(n_0,n_1)$. Writing the long exact sequence associated with 
\[\ol{F}_1^*\tilde{F}_0\hookrightarrow \ol{F}_1^*{F}_0\twoheadrightarrow \ol{F}^*_1\otimes\Tt\]
we find
\[
H^0(\ol{F}_1^*\tilde{F}_0)\hookrightarrow H^0(\ol{F}_1^*{F}_0)\longrightarrow H^0(\ol{F}_1^*\otimes\Tt)\longrightarrow H^1(\ol{F}_1^*\tilde{F}_0)\twoheadrightarrow H^1(\ol{F}_1^*{F}_0).
\]
By genericity, we may assume that $F_0$ and $F_1$ are stable. Thus, $H^0(\ol{F}_1^*{F}_0)\subset H^0({F}_1^*{F}_0)=0$, given that the stability of $\Ee$ implies that $\mu(F_0)<\mu(F_1)$. Hence, the sequence reads 
\[
H^0(\ol{F}_1^*\otimes\Tt)\hookrightarrow H^1(\ol{F}_1^*\tilde{F}_0)\twoheadrightarrow H^1(\ol{F}_1^*{F}_0).
\]
Let $e\in H^1(\lF_1^*\tF_0)$ be the extension which is the image of $f_1\in H^1(\tF_1^*\otimes\Tt)$. Then, by triviality of the pushforward extension,  $E\hookrightarrow F_0\oplus \ol{F}_1\stackrel{f_0+f_1}{\twoheadrightarrow}{\Tt}$. The induced  morphism $E\longrightarrow F_0$ is surjective, as otherwise $E\subset F_0'\oplus\ol{F}_1$, with $F_0'\subsetneq F_0$, which corresponds to an element in $H^0(\ol{F}_1^*\otimes \Tt')$ for some subsheaf $\Tt'\subsetneq \Tt$. Therefore $E\twoheadrightarrow F_0$ has kernel $F_1$ and we obtain that \eqref{eq:good_cd} holds.

Conversely, if $\Ee=(F_0\oplus F_1,\varphi)\in\Fix$ is of wobbly type $(n_0,n_1)$, let $\lim_{t\to 0}({E},\psi)=(F_0\oplus F_1,\varphi)$ for some $\psi$ such that $\lim_{t\to \infty}({E},\psi)=(F_0'\oplus F_1',\varphi')$, with $\rk(F_i')=n_i$. 
 Then, by Lemma \ref{prop:3.4_HH}, there are short exact sequences
\[
\xymatrix{
&F_0'\ar@{^(->}[d]&\\
F_1\ar@{^(->}[r]&E\ar@{->>}[d]\ar@{->>}[r]&F_0\\
&F_1'&
}.
\]
Then, since $\psi$ induces $\varphi'$   
 and $\psi|_{F_1}\equiv \varphi$  is injective for the generic fixed point $\Ee$, it follows that $F_1\cap F_0'=0$ and so $F_1\subset F_1'$. Given that ranks match $\Tt:=F_1'/F_1$ is torsion, and it must be $F_0/F_0':=\Tt$.  
Thus, 
there is a factorisation $\varphi:F_1\hookrightarrow F_1'=\ol{F}_1\stackrel{\psi}{\longrightarrow}F_0'K=\tilde{F}_0K\hookrightarrow F_0K$ where the saturation $\ol{F}_1$ and elementary modification $\tilde{F}_0$ are taken with respect to $\mathcal{T}$.
\end{proof}
\begin{lemma}
\label{lm:varphiperp}
Assume  $\varphi:F_1\longrightarrow F_0K$ is injective and $F_1$ is stable. Let $S:=F_0/{F}_1K^*$ be torsion free, where we identify $F_1K^*=\Im(\varphi)K^*$. 
Then $\varphi^\perp=\Im(H^1(S^*F_1)\longrightarrow  H^1(F_0^*F_1))$. Similarly, if $S$ is a torsion sheaf, then  $\varphi^\perp=\Im(H^0(S^*)\longrightarrow  H^1(F_0^*F_1))$
\end{lemma}
\begin{proof}
  By injectivity of $\varphi|_{F_1}$, $F_1\cap F_1K^*=0$, as $F_1K^*\subset\Ker(\varphi)$. When  $S$ is torsion free, there is a long exact sequence in cohomology 
    \begin{equation}
\label{eq:varphiperp}\cdots\longrightarrow H^1(S^*F_1)\longrightarrow H^1(F_0^*F_1)\twoheadrightarrow         \underbrace{H^1(\End(F_1)K)}_{\cong\CC}
    \end{equation}
    where the last isomorphism follows from stability of $F_1$. Similarly, if $S$ is torsion free or torsion, a similar exact sequence holds substituting $H^1(S^*F_1)$ by $H^0(S)$. Then, in both cases, $\varphi^\perp$ projects to either zero or the whole $H^1(\End(F_1)K)$. Given that $\varphi$ is injective, then $H^1(\End(F_1)K)=\varphi^*$, so the latter is not possible. Thus $\varphi^\perp\subset \Im(H^1(S^*F_1)\longrightarrow  H^1(F_0^*F_1))$ and similarly for the torsion case. Equality follows from hyperplanarity of $\varphi^\perp$.
\end{proof}

\begin{proposition}\label{prop:higher_order_wobbly}
    Let $n_0> n_1$. Let  $\Ee=(F_0\oplus F_1,\varphi)\in\Fix$ be general. In particular, we assume $\varphi$ to be injective, $S$ torsion free and $F_0$, $F_1$ stable. Suppose $\Ee$ is wobbly, and let $\lim_{t\longrightarrow 0}(E,\psi)=\Ee$. Then,  $E$ fits in an exact diagram
    \begin{equation}
        \label{eq:two_expressions_downwardflow}
        \xymatrix{F_1\ar@{^(->}[d]\ar@{=}[r]&F_1\ar@{^(->}[d]&\\
        F_1\oplus F_1K^*\ar@{^(->}[r]\ar@{->>}[d]& E\ar@{->>}[r]\ar@{->>}[d]&S\ar@{=}[d]\\
        F_1K^*\ar@{^(->}[r]& F_0\ar@{->>}[r]&S.}
    \end{equation}
    Moreover, if $\delta\leq n_0n_1(g-1)+1$, $E=F_0\oplus F_1$.
\end{proposition}
\begin{proof}

By Lemmas \ref{lm:perp} and \ref{lm:varphiperp}, $E\in\varphi^\perp=\Im(H^1(S^*F_1)\longrightarrow H^1(F_0^*F_1))$. In particular, $F_1K^*\hookrightarrow E$, as the extension class of $E$ maps to zero inside $H^1(F_1^*KF_1)$. This implies that $F_1\oplus F_1K^*\subset E$ is a vector subbundle. This together with Lemma \ref{prop:3.4_HH} prove the first statement. 

Note that $\varphi^\perp=0$ generically when $\delta\leq n_0n_1(g-1)+1$, as by Corollary \ref{cor:BN}, in this case $h^1(F_0^*F_1)=h^0(F_1^*F_0K)=1=h^1(\End(F_1)K)$, so that $\Im(H^1(S^*F_1)\longrightarrow H^1(F_0^*F_1))=0$. So by Lemma \ref{lm:perp}, in this case $E=F_0\oplus F_1$.
\end{proof}
\begin{proposition}\label{prop:wobbly_n0'n1'}  
   Let $n_0>n_1$, and assume that $\Fix$ is a component such that for the general point $H^0(S^*F_1K)=0$. Then, the set of points is of wobbly type $(n_0',n_1')\neq(n_0,n_1)$ is not dense.
\end{proposition}
\begin{proof}
    The proof is by dimensional count, by proving that a dense open set of points only intersects the set of fixed points of wobbly type $(n_0',n_1')$ in a lower dimensional subscheme. Let $\Ee$ be a generic fixed point satisfying the conditions and let $e\in H^1(F_0^*F_1)$ have underlying bundle $E$ admitting a nilpotent Higgs field $\psi$ with flow line joining $\Ee$ to a point $\Ee'=(F_0'\oplus F_1',\varphi')$ of type $(n_0',n_1')$ at infinity.
    Then $\psi^2=0$, as by \cite[Prop. 3.11]{HH} $F_0'\subset E_0=\Ker(\psi)$ and $\psi(E)\subset F_0'$. Thus we have an exact diagram
    \[
    \xymatrix{
&E_0\ar@{=}[r]\ar@{^(->}[d]&E_0\ar@{^(->}[d]\\
F_1\ar@{=}[d]\ar@{^(->}[r]&E \ar@{->>}[r]\ar@{->>}[d]&F_0\ar@{->>}[d]\\
F_1\ar@{^(->}[r]&E_1\ar@{->>}[r]&M.
}
 \]
Since $(n_0,n_1)\neq (n_0',n_1')$, $M$ cannot be a torsion sheaf by Proposition \ref{prop:wobbly_n0n1}. Now, if $M$ has torsion, since $E_1$ is torsion free, then we may saturate $F_1$ inside $E_1$ and $F_0$ inside $E_0$ to obtain 
    \[
    \xymatrix{
&E_0\ar@{^(->}[r]\ar@{^(->}[d]&\ol{E}_0\ar@{^(->}[d]\\
F_1\ar@{^(->}[d]\ar@{^(->}[r]&E \ar@{->>}[r]\ar@{->>}[d]&F_0\ar@{->>}[d]\\
\ol{F}_1\ar@{^(->}[r]&E_1\ar@{->>}[r]&M_{tf}.
}
    \]
    Now, by genericity, $\psi|_{F_1}=\varphi$ has no zeros. This means that $F_1=\ol{F}_1$, $E_0=\ol{E}_0$, namely, $M=M_{tf}$. 
    
Then, if $H^0(M^*{E}_0K)\neq 0$, then $F_0$ admits an order two nilpotent Higgs field, namely, it is wobbly, so non generic by Lemma \ref{lm:dominance}. So we may assume $H^0(M^*{E}_0K)=0$. Since $F_0$ is assumed to be general, then it must be  
$d_{M}n_{E_0}-n_{M}n_{E_0}\geq n_{E_0}n_{M}(g-1)$,
by \cite[Theorems 0.1\&0.2]{RT}. Then
    $\psi\in H^0(E_1^*E_0K)\subset H^0(F_1^*E_0K)$. Namely, we may identify $\psi=\varphi$. Since $E_0\oplus F_1\subset F_0\oplus F_1$, it follows that $\Ee$ is an extension of the trivial Higgs bundle $(M,0)$ by $(E_0\oplus F_1,\varphi)$. These are parametrised by the first hypercohomology of the complex 
    \[
   C_\bullet\,:\, M^*(E_0\oplus F_1)\longrightarrow  M^*(E_0\oplus F_1)K\qquad s\mapsto \varphi\circ s.
    \]
    Using the exact sequence
    \[
     M^*E_0\hookrightarrow M^*(E_0\oplus F_1)\longrightarrow  M^*(E_0\oplus F_1)K\twoheadrightarrow  M^*E_0K\oplus \overbrace{M^*F_1K/M^*F_1}^{\Tt:=}
    \]
that we may reinterpret in terms of the short exact sequence of complexes 
\[
\overbrace{\left[M^*E_0\longrightarrow 0\right]}^{C'_\bullet}\hookrightarrow 
 {\left[M^*(E_0\oplus F_1)\longrightarrow  M^*(E_0\oplus F_1)K\right]}\twoheadrightarrow \overbrace{\left[M^* F_1\longrightarrow  M^*(E_0\oplus F_1)K\right]}^{C_\bullet^{''}}\]
 together with a quasi isomorphism
 \begin{equation}
     \label{eq:qi}
 \left[M^* F_1\longrightarrow  M^*(E_0\oplus F_1)K\right]\cong_{qi}\left[0\longrightarrow   M^*E_0K\oplus \Tt\right]
 \end{equation}
 we see that the hypercohomology 
   fits in a short exact sequence
   \[
   H^1(M^*E_0)\hookrightarrow \HH^1(C_\bullet)\twoheadrightarrow H^0(M^*E_0K\oplus \Tt)=H^0(\Tt)
   \]
   where we have used \eqref{eq:qi} to deduce $\HH^0({C_\bullet}^{''})=0$, as well as $\HH^2(C_\bullet')=0$, and $H^0(M^*E_0K)=0$. Letting $\delta_0:=n_1\deg(E_0)-d_1n_{E_0}+2n_{E_0}n_1(g-1)$, the same dimensional computation as the one leading to \eqref{eq:dim_Fix} shows that the family parametrised by these objects will be of dimension at most
   \[
   -1+h^0(M^*E_0K\oplus M^*F_1K/M^*F_1)+H^1(M^*F_1)+h^0(E_0F_1^*K)+(n_1^2+n_M^2+n_{E_0}^2)(g-1)+3\]
   \[=2+n_Mn_1(2g-2)-\chi(M^*F_1)+\min\{\chi(E_0F_1^*K),1\}+(n_1^2+n_0-2n_Mn_{E_0})(g-1)\]
   \begin{equation}
       \label{eq:some_points}
   =\left\{\begin{array}{ll}
    2+\delta-\delta_0++(n_1^2+n_0-2n_Mn_{E_0})(g-1)    &  \textrm{ if }\delta_0>n_1n_{E_0}(g-1)\\
     3+\delta-\delta_0++(n_1^2+n_0-2n_Mn_{E_0})(g-1)   & \textrm{ if }\delta_0<n_1n_{E_0}(g-1).
   \end{array}\right.\end{equation}
   In the above we have used $\delta-\delta_0=n_Mn_1(2g-2)-\chi(M^*F_1)$. Since, $F_0$ is stable, then $\delta_0<\delta$, and so the dimension in \eqref{eq:some_points} is strictly smaller than \eqref{eq:dim_Fix}. 
\end{proof}

\section{Wobbliness of length two components}\label{sec:wobbly_comps}
\subsection{$\GL_n(\CC)$-wobbliness of length two components}
In this section we prove that all fixed points of type $(n_0,n_1)\neq (2,1), (1,1)$ are wobbly.
\begin{theorem}\label{thm:wobbly_comps}
Let $\delta\in\IC$. Let $n=n_0+n_1>3,\, \ol{n}=(n_0,n_1)$ and assume $n_0\geq n_1$. Then, any fixed point component $\Fix$ is wobbly. Moreover, given a general point $\Ee=(F_0\oplus F_1,\varphi)$, then 
\begin{enumerate}
    \item\label{it:wobbly_delta<3n1(n0-n1)} If $\delta<3n_1(n_0-n_1)(g-1)$ then $\Ee$ is wobbly, realised by a Higgs field of the form $\psi=\begin{pmatrix}
        0&\varphi\\\beta&0
    \end{pmatrix}$
    for some $\beta\in H^0(S^*F_1K)$ with limit at infinity of order three equal to  $(F_1K^*\oplus F_1\oplus S,\psi)$ if the latter were semistable. 
    \item\label{it:wobbly_delta>3n1(n0-n1)} If $\delta\geq 3n_1(n_0-n_1)(g-1)$ then $\Ee$ with flow realised by a Higgs field of order four of the form $\begin{pmatrix}
        0&\varphi\\\beta&0
    \end{pmatrix}$
    for some $\beta\in H^0(F_0^*F_1K)$ not factoring through $S:=F_0/F_1K^*$.  
\end{enumerate}
\end{theorem}
\begin{remark}\label{rk:cases_thm_wobbly}
   Since $\delta\in\IC$, it follows that  when $n_0=n_1$, only Theorem \ref{thm:wobbly_comps} \eqref{it:wobbly_delta>3n1(n0-n1)} applies. Indeed, $\delta<0=3n_1(n_0-n_1)(g-1)$ does not correspond to any fixed point component. Similarly, if $n_1=1$ $n_0\neq 2$, then $3n_1(n_0-n_1)\geq 2n_0n_1$,  Theorem \ref{thm:wobbly_comps} \eqref{it:wobbly_delta>3n1(n0-n1)} is empty and Theorem \ref{thm:wobbly_comps} \eqref{it:wobbly_delta<3n1(n0-n1)} covers all cases. The argument breaks for rank $3$, and is treated in joint work with Pauly \cite{PPrk3}.
\end{remark}
Before proving the statement, let us look at the key example:
\begin{example}\label{ex:key_example}
    Consider a Higgs bundle of the form $(FK^*\oplus F,1)$ with the tautological Higgs bundle $1:F\longrightarrow FK^*\otimes K$. Then, if $\rk(F)\neq 1$, for a general $F$, there exists a $K^2$-twisted nilpotent Higgs bundle \cite[Lemma 5.2]{W=S}. Thus, if $\beta\in H^0(\End(F)K^2)$ is nilpotent, it follows that the Higgs field $\begin{pmatrix}0&1\\
    \beta&0
    \end{pmatrix}$ is a nilpotent Higgs field on $FK^*\oplus F$. Since $\left(FK^*\oplus F,t\begin{pmatrix}0&1\\
    \beta&0
    \end{pmatrix}\right)$ is gauge equivalent to $\left(FK^*\oplus F,\begin{pmatrix}0&1\\
    t\beta&0
    \end{pmatrix}\right)$, it follows that $FK^*\oplus F$ is wobbly.
\end{example}
Now, for a general fixed point $(F_0\oplus F_1,\varphi)$, $(F_1K^*\oplus F_1,1)\subset (F_0\oplus F_1,\varphi)$ is a Higgs subbundle by Proposition \ref{prop:higher_order_wobbly}. The remainder of this section is devoted to proving that for $n_1>1\delta$ large enough, the Higgs field from Example \ref{ex:key_example} lifts to $F_0\oplus F_1$. 
\begin{proof}[Proof of Theorem \ref{thm:wobbly_comps}]
By Proposition \ref{prop:vs_are_dense}, it is enough to prove that the general point is wobbly for the given range.

First note that for $\delta< 3n_1(n_0-n_1)(g-1)$, it follows that there exists a non-zero
\[
\beta:F_0\longrightarrow F_1K\textrm{ s.t. }\beta\circ\varphi=0.
\]
Indeed, let $\tilde{S}$ be the torsion free subsheaf of $F_0/F_1K^*$. Then 
\[
\deg(\tilde{S})\leq d_0-d_1+2n_1(g-1)
\]
so that
\[
\deg(\tilde{S}^*F_1K)\geq (n_0-n_1)d_1+2(n_0-n_1)n_1(g-1)-n_1(d_0-d_1+2n_1(g-1))\]\[
=n_0d_1-n_1d_0+2n_1(n_0-2n_1)(g-1)=\]\[=-\delta+4n_1(n_0-n_1)(g-1).
\]
So by Riemann--Roch,
\[
\chi(S^*F_1K)=-\delta+4n_1(n_0-n_1)(g-1)-(n_0-n_1)n_1(g-1)>0
\]
for all pairs $(S,F_1)$ whenever $\delta< 3n_1(n_0-n_1)(g-1)$, so $ 
h^0(S^*F_1K)>0$. This proves wobbliness of all such components. The statement about limits follows from Proposition \ref{prop:3.4_HH}.

Assume now that $\delta\geq 3n_1(n_0-n_1)(g-1)$ (so that $n_1\neq 1$ or $n_0=2$ by Remark \ref{rk:cases_thm_wobbly}). 
 By Lemma \ref{lm:dominance}, we may assume  $F_1$ to be generic. By \cite[Thm. 0.1]{RT} (see also \cite[Thm. 4.4]{Hirsch}, for every ordered partition $(n_M,n_N)$ of $n_1$, $F_1$ is an extension of the form
\[
M\hookrightarrow F_1\twoheadrightarrow N
\]
with $M,\ N$ stable bundles of ranks $\rk(M)=:n_M$, $\rk(N)=:n_N$, and degrees $\deg(M)=:d_M$, $\deg(N)=:d_N$ satisfying 
\begin{equation}
    \label{eq:dNrM}
    0<\mu(N)-\mu(M)= (g-1)+\frac{\ell}{n_Mn_N}
\end{equation}
where $\ell=0,\dots, n_1-1$.
Then, we have the following exact diagram:
\begin{equation}
    \label{eq:diag_MN}
    \xymatrix{
MK^*\ar@{=}[r]\ar@{^(->}[d]&MK^*\ar@{^(->}[d]&\\
F_1K^*\ar@{^(->}[r]\ar@{->>}[d]&F_0\ar@{->>}[r]\ar@{->>}[d]&\ar@{=}[d]S\\
NK^*\ar@{^(->}[r]&\tilde{N}K^*\ar@{->>}[r]&S.
}
\end{equation} 
 We want to prove that $h^0((\tilde{N}K^*)^*MK)>0$, as then, the composition $F_0\twoheadrightarrow \tilde{N}K^*\stackrel{\beta}{\longrightarrow }MK\hookrightarrow F_1K$ satisfies that the Higgs field $\begin{pmatrix}
    0&\varphi\\
    \beta&0
\end{pmatrix}$ is nilpotent and fulfils the conditions of Lemma \ref{prop:3.4_HH}. 

Note that 
\[
h^0(\tilde{N}^*MK^2)\geq n_{\tilde{N}}n_M(\mu(M)-\mu(\tilde{N})+3g-3),\]
so it is enough to show that
\begin{equation}
    \label{eq:key_ineq}
\mu(\tilde{N})-\mu(M)<3g-3.
\end{equation}

Without loss of generality, we may assume that $d<n$ (as up to $\CC^\times$-equivariant isomorphism, the moduli space satisfies this), so that $d_1>0$, namely, $d_1\geq 1$, and $d_0<1$, namely, $d_0\leq 0$.

Without loss of generality, we may assume that $d<n$ (as up to $\CC^\times$-equivariant isomorphism, the moduli space satisfies this), so that $d_1>0$, namely, $d_1\geq 1$, and $d_0<1$, namely, $d_0\leq 0$. Let 
\[ 
n_M=1\Rightarrow n_N=n_1-1,\, n_{\tN}=n_0-1.
\]

We may rewrite \eqref{eq:dNrM} as
\[
d_{{N}}n_M-d_Mn_{{N}}=d_1-d_Mn_1=(g-1)(n_1-1)+\ell,\]
or, equivalently,
\begin{equation}
    \label{eq:dNrM2}
    d_M=\frac{1}{n_1}\left(d_1-\ell-(g-1)(n_1-1)\right).
\end{equation}
Similarly, \eqref{eq:key_ineq} may be rewritten as
\begin{equation}
    \label{eq:key_ineq2}
d_{\tilde{N}}n_M-d_Mn_{\tilde{N}}=d_0-d_Mn_0<3(n_0-1)(g-1)\iff d_0<d_Mn_0+3(n_0-1)(g-1).
\end{equation}
so, to prove \eqref{eq:key_ineq2}, it suffices to show that
\[d_Mn_0+3(n_0-1)(g-1)>0.\]
Plugging \eqref{eq:dNrM2} into the left hand side of the previous inequality and reorganising the terms, we find
\[
d_Mn_0+3(n_0-1)(g-1)=(g-1)\left(3n_0-3-\frac{n_0(n_1-1)}{n_1}\right)+\frac{n_0}{n_1}(d_1-\ell)
\]
\[
\stackrel{\tiny d_1>0}{\geq} (g-1)\left(2n_0-3+\frac{n_0}{n_1}\right)+\frac{n_0(1-\ell)}{n_1}
\]
\begin{equation}
    \label{eq:final_ineq}
    \stackrel{\ell\leq n_1-1}{\geq} (g-1)\left(2n_0-3+\frac{n_0}{n_1}\right)+\frac{n_0(2-n_1)}{n_1}
\end{equation}
Since $n_0\geq 2$, $2n_0-3\geq n_0$ for $n_0>2$, and $2n_0-3=1$ if $n_0=2$, in which case $n_1=2$ (as $n_1>1$, cf Remark \ref{rk:cases_thm_wobbly}). In the latter case, we have that the right hand side of \eqref{eq:final_ineq} equals
\[
2(g-1)>0\geq d_0
\]
and we are done. If $n_0\geq 2$, so that $2n_0-3\geq n_0$, then the right hand side of \eqref{eq:final_ineq} is bounded by
\[
\eqref{eq:final_ineq}\geq (g-1)\left(n_0+\frac{n_0}{n_1}\right)+\frac{n_0(2-n_1)}{n_1}=(g-1)\frac{n_0}{n_1}+\frac{2n_0}{n_1}>0\geq d_0
\]
which finishes the proof.
\end{proof}
\begin{remark}
    \label{rk:segre_strata}
   Note that a field $\beta\in H^0(\tN^*MK^2)$ extends the field $\beta|_N\in H^0(N^*MK^2)\setminus 0$ constructed in Example \ref{ex:key_example}, which is a $K^2$-twisted Higgs field associated to the Segre stratum of $F_1$, as .
\end{remark}
\begin{remark}
 Theorem \ref{thm:wobbly_comps} explains some phenomena observed for particular cases. For instance, in rank four, all components of type $(3,1)$ were known to be wobbly \cite[Remark 5.13]{HH}, which is the first fully wobbly case covered by the criterion in Theorem \ref{thm:wobbly_comps}, corresponding to $n_1=1$.
\end{remark}
\begin{remark}\label{rk:limit_infty_delta_big}
       The proof of Theorem \ref{thm:wobbly_comps} provides a candidate for the limit at infinity also in the case $\delta\geq 3n_1(n_0-n_1)(g-1)$. 
       With the notation of Theorem \ref{thm:wobbly_comps}, let $\beta\in H^0((\tilde{N}K^*)^*MK)$. Then, if the point \[\Ee':=\left({M}K^*\oplus M\oplus \tilde{N}K^*\oplus N,\psi:=\begin{pmatrix}
    0&1&0&0\\
    0&0&\beta&0\\
    0&0&0&1\\
    0&0&0&0
\end{pmatrix}\right)\]
is semistable, then by Lemma \ref{prop:3.4_HH}
\[\Ee'=\lim_{t\to \infty}t\cdot\left(F_0\oplus F_1,\begin{pmatrix}
   0&\varphi\\
\beta&0
\end{pmatrix}\right).\]
In this case, also by Lemma \ref{prop:3.4_HH}, $\Ee$ is of wobbly type $(n_0-1,n_1-1,1,1)$. 
  \end{remark} 
\begin{remark}\label{rk:wobbly_type}
As a remark to Theorem \ref{thm:wobbly_comps}, we see that the wobbly type is not uniquely defined: when $n_1>1$ and $\delta<3n_1(n_0-n_1)(g-1)$, there are at least two order two Higgs fields on $\Ee$. Similarly, by Theorem \ref{thm:wobbly_comps} and Proposition \ref{prop:wobbly_n0n1}, points of wobbly type $(n_0,n_1)$ characterised in Proposition \ref{prop:wobbly_n0n1}, are also of wobbly type $(m_1,\dots,m_s)$ for $s\geq 3$.
\end{remark}


\subsection{$\U(n_0,n_1)$-wobbliness of length two components.}
This section gathers some observations resulting in Corollary \ref{cor:wobblyiffrealwobbly}.
\begin{definition}\label{def:real_vs}
Let $G_\mathbb{R}<\GL(n,\CC)$ be a real form. A fixed point $\Ee\in\M_X(G_\RR)$ is $G_\mathbb{R}$-very stable if there exists no nilpotent $G_\mathbb{R}$-Higgs bundle $(E,\psi)$ such that $\lim_{t\to 0}(E,\psi)=\Ee$. Otherwise it is called $G_\mathbb{R}$-wobbly. 

Similarly,  a fixed point component $\Fix\subset\M_X(G_\RR)$ is called $G_\mathbb{R}$-very stable, if the generic point is $G_\mathbb{R}$-very stable, and $G_\mathbb{R}$-wobbly otherwise.
\end{definition}
Note that $\Fix\subset\M_X(\U(n_0,n_1))$. Since $\M_X(\U(n_0,n_1))\cong \M_X(\U(n_1,n_0))$, when speaking of real Higgs bundles we assume $n_0\geq n_1$. 

More generally, for an arbitrary complex reductive Lie group $G$, let $\M_X(G)$ be the moduli space of $G$-Higgs bundles, and $h_G: \M_X(G)\longrightarrow B_G$ the corresponding Hitchin map. Then,
\[h_G^{-1}(0)\subset \bigcup_{G_\RR<G Hodge}\M_X(G_\RR)\]
where the union runs over all real forms $G_\RR<G$ of Hodge type \cite{Simpson}.
\begin{remark}\label{rk:wobblyiffgricwobblyreal}
    Proposition \ref{prop:vs_are_dense} is still valid if we substitute wobbliness by $G_\RR$ wobbliness. Indeed, the downward flow intersected with $\M_X(G_\RR)$ and the nilpotent cone is still locally of finite type by closedness of $\M_X(G_\RR)\subset \M_X(n,d)$. Nevertheless, not all the Higgs bundles realising wobbliness realise $G_\mathbb{R}$-wobbliness. As an example, for $G_\RR=\U(n_0,n_1)$, when $F_0/\varphi(F_1)K^*$ has torsion, there are flow lines  outside of $\M_X(\U(n_0,n_1))$ flowing down to $(F_0\oplus F_1,\varphi)$.
\end{remark}

Clearly, if a component is very stable, then it is  $G_\mathbb{R}$-very stable. As a corollary to Theorem \ref{thm:wobbly_comps} we find that for length two fixed points this is an equivalence. 
\begin{corollary}\label{cor:wobblyiffrealwobbly}
The fixed point component $\Fix$ is wobbly  if and only if it is $\U(n_0,n_1)$-wobbly.
\end{corollary}
\begin{proof}
By Theorem \ref{thm:wobbly_comps}, a component is wobbly if and only if it is flown down into by flows of the form $(F_0\oplus F_1,\psi)$ where
\[
\psi=
\begin{pmatrix}
0&\beta\\
\varphi&0\end{pmatrix}.
\]
This is by definition an $U(n_0,n_1)$-Higgs bundle.
\end{proof}
Corollary \ref{cor:wobblyiffrealwobbly} yields a simpler way to check whether a component is or not very stable, as finding (or refuting the existence of) downward flows inside $\M_X(\U(n_0,n_1))$is noticeably easier than doing so inside $\M_X(n,d)$. Another interesting consequence of Corollary \ref{cor:wobblyiffrealwobbly} is the following:
\begin{corollary}
    Let $n_0\geq n_1$, and denote by $\M^{\tau}(\U(n_0,n_1))$ the connected component with Toledo invariant $\tau$ such that $0<|\tau|\leq 2\min\{p,q\}(g-1)$, with equality if and only if $p=q$. If $n_1\neq 1$, then $\M^{max}(\U(n_0,n_1))$ hits the nilpotent cone in different fixed point components. 
\end{corollary}
\begin{proof}
    By Theorem \ref{thm:description_fixed_points}, when $n_0=n_1$, the order two fixed point component $\Fix$ for is non empty, and wobbly by Theorem \ref{thm:wobbly_comps}, so contained in the $\U(n_0,n_1)$-wobbly locus by Remark \ref{rk:wobblyiffgricwobblyreal}. By Corollary \ref{cor:wobblyiffrealwobbly}, since $n_1\neq 1$, there is a downward flow into it inside $\M_X(\U(n_0,n_1))$. We observe that the Toledo invariant inside the flow line is the same. By connectedness, it can only flow up within the maximal component of the Toledo invariant. 

    For $n_0>n_1>1$, $\M^{max}(\U(n_0,n_1))\cong \M^{max}(\U(n_1,n_1))\times \N(n_0-n_1,d')$ by \cite[Thm. 3.32]{BGGP}. So we may apply the preceding result to $\M^{max}(\U(n_1,n_1))$ using $n_1\neq 1$. 
\end{proof}
\section{Virtual equivariant multiplicities of wobbly components}\label{sec:multiplicities}

In what follows, we analyse the failure of polynomiality of equivariant multiplicities.

\begin{lemma}\label{lm:weights}
Let $\Ee=(F_0\oplus F_1,\varphi)$ be a smooth fixed point. Then the $\CC^\times$-weight subsheaves of  on $\End(E)$ are as described below:
   \begin{table}[h]
       \centering
        \begin{tabular}{|c|c|c|c|c|}
    \hline
        Subsheaf & $F_1^*F_0$& $F_0^*F_1$&$\End(F_0)$&$\End(F_1)$ \\\hline
         Weight & 1&-1&0&0\\\hline
    \end{tabular}
     \end{table}
    \end{lemma}
\begin{proof}
The weight $w_i$ of the action on $F_i$ can be computed by imposing $n_0w_0+n_1w_1=0$, and $w_0-w_1=1$. This implies that the gauge transformation 
\[
Ad\begin{pmatrix}
t^{\frac{n_1}{n}}Id_{F_0}&\\
&t^{\frac{-n_0}{n}}Id_{F_1}
\end{pmatrix}
\]
takes $\Ee$ to $t\Ee$. Thus, the weight on $F_0$ is $w_0=\frac{n_1}{n}$ and the weight on $F_1$ is $w_1=\frac{-n_0}{n}$.
One easily computes from that the weights on the subsheaves $F_iF_j^*$. 
\end{proof}
\begin{proposition}\label{prop:weight_spaces}
Let $\Ee=(F_0\oplus F_1,\varphi)=:(E,\varphi)$ be general; in particular, assume it is smooth and $\varphi$ is injective.  Then, the $\CC^\times$-module $T_\Ee^+$ consists of the following weight spaces:
\[
(T_\Ee^+)_1\cong \frac{H^0(F_1^*F_0K)}{H^0(\End(F_1))\oplus\End(F_0))}\oplus  \left(\frac{H^0(\End(F_1)K)\oplus\End(F_0)K)}{H^0(F_0^*F_1)}\right)^*  
\]
\[
(T_\Ee^+)_2\cong H^1(F_0F_1^*)
\]
where $S:=F_0/F_1K^*$ is locally free if $n_0>n_1$ and $\mathcal{T}:=F_0/F_1K^*$ is torsion otherwise. In particular
\[
\dim (T_\Ee^+)_1=\left\{\begin{array}{ll}
\delta+1+n_0^2(g-1)   &  n_0=n_1;\\
  (n_0^2+n_1^2-n_0n_1)(g-1)+1+\delta& n_0>n_1; 
\end{array}
\right.
\]
\[
\dim (T_\Ee^+)_2=-\delta+3n_0n_1(g-1).
\]
\end{proposition}
\begin{proof}
Recall that $T_\Ee=\HH^1(C_\bullet)$ for the complex
\[
C_\bullet:\End(F_0\oplus F_1)\stackrel{\textrm{ad}(\varphi)}{\longrightarrow} \End(F_0\oplus F_1)K.
\]
Hence, using the long exact sequence in hypercohomology defining $\HH^1$, we have a short exact sequence 
\begin{equation}\label{eq:SES_defining_HH1}
    \xymatrix{
\Coker\left(H^0(\End(E))\stackrel{h^0(\textrm{ad}(\varphi))}{\rightarrow}H^0(\End(E)\otimes K)\right)\ar@{^(->}[d]\\ T_\Ee\ar@{->>}[d]\\ \Ker\left(H^1(\End(E))\stackrel{h^1(\textrm{ad}(\varphi))}{\rightarrow}H^1(\End(E)\otimes K)\right).}
\end{equation}
Now, we have the following exact sequences:
\begin{eqnarray}\label{eq:exact_SESs_n0>n1}
F_0F_1^*\hookrightarrow F_0F_1^*\stackrel{ad[\varphi]}{\twoheadrightarrow} 0;
\\\nonumber F_0^*F_1\stackrel{ad[\varphi]}{\hookrightarrow} \End(F_0)K\oplus \End(F_1)K\twoheadrightarrow\left\{\begin{array}{ll}
 \End(F_0)K{\oplus}F_1\otimes\mathcal{T}    &  n_0=n_1;\\
 F_0^*SK \oplus\End(F_1)K   & n_0>n_1;
\end{array}\right.
\\\nonumber
\left.\begin{array}{ll}0  &  n_0=n_1;\\
 S^*F_0    & n_0>n_1;
\end{array}\right\}\hookrightarrow \End(F_0) 
  \stackrel{ad[\varphi]}{\longrightarrow } F_0F_1^*K\twoheadrightarrow\left\{\begin{array}{ll}
 F_0\otimes\mathcal{T}   &  n_0=n_1;\\
 0    & n_0>n_1.
\end{array}\right.\\\nonumber \End(F_1)\stackrel{ad[\varphi]}{\hookrightarrow}  F_1^*F_0K\twoheadrightarrow\left\{\begin{array}{ll}
 F_0K\otimes\mathcal{T}  &  n_0=n_1\\
 F_1^*\otimes S\otimes K   & n_0>n_1.
\end{array}\right.
\end{eqnarray}
Hence, from  \eqref{eq:SES_defining_HH1},  \eqref{eq:exact_SESs_n0>n1},  
 and \cite[Prop. 4.1]{Gothen} one sees that the  positive weight contributions to $\HH^1(C_\bullet)$  come from the weight one space
\[
  \frac{H^0(F_1^*F_0K)}{H^0(\End(F_1))\oplus\End(F_0))}\oplus  \left(\frac{H^0(\End(F_1)K)\oplus\End(F_0)K)}{H^0(F_0^*F_1)}\right)^* 
\]
and the weight two space $H^1(F_0F_1^*)$. The latter has dimension
\[
\dim(T_\Ee^+)_2=h^1(F_0F_1^*)=-\delta+3n_0n_1(g-1)
\]
by stability of $\Ee$ and unstability of $F_0\oplus F_1$. 

To compute the exact dimensions of the weight one spaces, by Corollary \ref{cor:BN} we have that generically
\[
h^0(F_0^*F_1)_{gen}=\left\{\begin{array}{ll}
     0& \textrm{ if } \delta>n_0n_1(g-1),\\
     1+n_0n_1(g-1)-\delta& \textrm{ if } \delta\leq n_0n_1(g-1),
\end{array}\right.
\]
so 
\[
\dim\left(\frac{H^0(\End(F_1)K)\oplus\End(F_0)K)}{\tiny H^0(F_0^*F_1)}\right)^* =\]\[=\left\{\begin{array}{ll}
     2+(n_0^2+n_1^2)(g-1)& \delta>n_0n_1(g-1),\\
     1+(n_0^2+n_1^2-n_0n_1)(g-1)+\delta&   \textrm{otherwise}.
\end{array}\right.
\]
In the above we have used $H^1(\End(F_i)K)=1$ (by stability of $F_i$ when $\delta\neq n_1(n_0-n_1)(g-1)$ and for the general pair $(F_0,F_1)$ otherwise by Theorem \ref{thm:description_fixed_points}).

For $n_0>n_1$ we have that $\Im(\End(F_1)\oplus\End(F_0)\longrightarrow F_1^*F_0K)=\Im(\End(F_0)\longrightarrow F_1^*F_0K)$ by \eqref{eq:exact_SESs_n0>n1}. For $n_0=n_1$, we have that $\Im(H^0(\End(F_1)\oplus\End(F_0)\longrightarrow H^0(F_1^*F_0K))=\CC$ by stability of $F_1$ and $S$. So
\[
\dim \frac{H^0(F_1^*F_0K)}{H^0(\End(F_0)\oplus\End(F_1))}=\left\{\begin{array}{ll}
     \delta-n_0n_1(g-1)-1& \textrm{ if } \delta>n_0n_1(g-1);\\
     0& \textrm{ otherwise }.
\end{array}\right.
\]
Thus, we have
\[
\dim (T_\Ee^+)_1=1+(n_0^2+n_1^2-n_0n_1)(g-1)+\delta
\]
which finishes the proof.
\end{proof}
\begin{proposition}\label{prop:equivar_mult}
Let $\Ee$ be a fixed point in a component of type $(n_0,n_1)$ with invariant $\delta$. Then 
\[
m_\Ee(t)=[2]_t^{e(\delta,n_0,n_1)}{p(t)}
\]
where 
$[m]_t:=\frac{1-t^m}{1-t}$ is the $m$-th quantum integer,
\[
e(\delta,n_0,n_1)=(g-1)\left(-3n_0n_1+2\left\lfloor\frac{n}{2}\right\rfloor^2+\left\lfloor\frac{n}{2}\right\rfloor\right)+\delta,
\]
and
\[
p(t)=\prod_{2\not|m\leq n}[m]_t^{(2m-1)(g-1)}\prod_{2|m\leq n}\left(\frac{[m]_t}{[2]_t]}\right)^{(2m-1)(g-1)}
\]
 is a product of cyclotomic polynomials coprime with $(1-t)(1+t)$.
\end{proposition}
\begin{proof}
By Proposition \ref{prop:weight_spaces}, and the fact that $(1-t^2)=(1+t)(1-t)$, we have
\[
m_\Ee(t)=\frac{(1-t)^g\prod_{i=2}^n(1-t^i)^{(2i-1)(g-1)}}{(1-t)^{\dim(T_\Ee^+)_1-\delta+3n_0n_1(g-1)}(1+t)^{-\delta+3n_0n_1(g-1)}}.
\]
Now, by the properties of cyclotomic polynomials, we have that $1-t^i=(1-t)(1+\dots +t^{i-1})$ and 
\[
 \left((1+t),(1-t^{2^km})\right)=\left\{\begin{array}{ll}
    1+t   & \textrm{ if } k>0,\ (m,2)=1 \\
    1  & \textrm{ if } k=0,\ (m,2)=1. 
 \end{array}\right.
\]
Thus, we obtain that
\[
m_\Ee(t)=(1-t)^{e_1}(1+t)^{e_2}p(t),
\]
where $p(t)$ is as stated, and
\[
e_1=g+\sum_{i=2}^n(2i-1)(g-1)-\dim(T_\Ee^+)_1-\dim(T_\Ee^+)_2\]
\[=1+n^2(g-1)-\dim(T_\Ee^+)=0
\]
where we have used that downward flows have Lagrangian fibers. On the other hand 
\[
e_2=\sum_{2\leq 2i\leq n}(4i-1)(g-1)=\left(2\left\lfloor\frac{n}{2}\right\rfloor^2+\left\lfloor\frac{n}{2}\right\rfloor\right)(g-1)+\delta-3n_0n_1(g-1).
\]
\end{proof}
Next we analyse the detection range of wobbliness given by polynomiality of equivariant multiplicities.

\begin{corollary}
    \label{thm:mults} Let $\Ee\in\Fix$ be a smooth fixed point. Then $m_\Ee(t)$ is not a polynomial if and only if
\[
\delta<\left(3n_0n_1-\left(2\left\lfloor\frac{n}{2}\right\rfloor^2+\left\lfloor\frac{n}{2}\right\rfloor\right)\right)(g-1).
\]
This range intersects the range $\IC$ non emptily for ordered partitions $(n_0,n_1)$ such that
\begin{equation}
    \label{eq:detection_range}
    n_1>n-\sqrt{n^2-2\left\lfloor\frac{n}{2}\right\rfloor^2-\left\lfloor\frac{n}{2}\right\rfloor}.
\end{equation}
\end{corollary}
\begin{proof}
The first statement follows from Proposition \ref{prop:equivar_mult}.

To see the second assertion, consider the function of $n_1$
\[
f_n(n_1)={3(n-n_1)n_1}-{\left(2\left\lfloor\frac{n}{2}\right\rfloor^2+\left\lfloor\frac{n}{2}\right\rfloor\right)}
.
\]
Then $e_2(n_0,n_1,\delta)=\delta-f_{n_0+n_1}(n_1)(g-1)$. The range $e<0$ is non empty if and only if
\[
f_n(n_1)>\left\{\begin{array}{ll}
n_1(n_0-n_1)     & \textrm{ if } n_0=n_1 \\
 n_1(n_0-n_1)+1    & \textrm{ if } n_0\neq n_1 
\end{array}\right..
\]
Now, let $p_n(n_1):=-n_1^2+2nn_1-2\left\lfloor\frac{n}{2}\right\rfloor^2-\left\lfloor\frac{n}{2}\right\rfloor$. Then, the above is equivalent to 
\[
p(n_1)> 0 \textrm{ if } n_0=n_1 \qquad p(n_1)>1 \textrm{ if }  n_0\neq n_1 
\]
We find that the only root of $n_1^2-2nn_1+2\left\lfloor\frac{n}{2}\right\rfloor^2+\left\lfloor\frac{n}{2}\right\rfloor$ smaller than $n$ is $\lambda_c:=n-\sqrt{n^2-2\left\lfloor\frac{n}{2}\right\rfloor^2-\left\lfloor\frac{n}{2}\right\rfloor}$, which yields the result by noticing that $\lambda_c<\left\lfloor\frac{n}{2}\right\rfloor$ and for $n>2$ the partition $(\left\lceil\frac{n}{2}\right\rceil,\left\lfloor\frac{n}{2}\right\rfloor)$, $f_n(\left\lfloor\frac{n}{2}\right\rfloor)>\delta_{min}$.
\end{proof}
\begin{remark}\label{rk:rk_3_full}
In rank $3$, the only fixed point components which are not of type $(1,1,1)$ are of type $(2,1)$ or $(1,2)$. It is proven in \cite{PPrk3} that, in this case, wobbliness is totally determined by non polynomiality of the virtual equivariant multiplicities. 
 This was conjectured in \cite[Remark 5.12]{HH}.  Similarly, Theorem \ref{thm:mults} implies that for partitions of the form $(n-1,1)$, $n>3$ the range \eqref{eq:detection_range} does not intersect $\IC$, which is particular explains why in rank $4$, components of type $(3,1)$ have always polynomial equivariant multiplicities (as noticed in \cite[Remark 5.13]{HH}).
 \end{remark}


\section{Euler pairings}\label{sec:Euler_pairings}
This final section computes a symmetrised equivariant Euler pairing between certain $\CC^\times$-equivariant branes. These branes are, on the one hand, the $BAA$-brane given by the structure sheaf on downward flows to type $(n_0,n_1)$ fixed points, and on the other, the $BBB$ branes proposed in \cite{HH} to be the dual to downward flows to very stable Higgs bundles of type $(1,\dots,1)$. 

\begin{definition}\label{def:equivar_euler_pairings}Given $\Ss_1$, $\Ss_2$ $\CC^\times$-equivariant sheaves on a semiprojective variety $M$, their equivariant Euler pairing is
\begin{equation}
    \label{eq:equivar_euler_pairings}
\chi_{\CC^\times}(M,\Ss_1,\Ss_2):=\sum_{i,j}(-1)^i H^i(M,RHom(\Ss_1,\Ss_2))^jt^{-j},
\end{equation}
where $H^i(M,RHom(\Ss_1,\Ss_2))^j\subset H^i(M,RHom(\Ss_1,\Ss_2))$ is the subspace under which $\CC^\times$ acts with weight $j$.
\end{definition}
Let
\begin{equation}
    \label{eq:Ff}
    \Ff:=(L\oplus L(D_1)K^{-1}\oplus\dots LK^{-n+1}(D_1+\dots+D_{n-1}),\bigoplus \phi_i)
\end{equation}
be a type $(1,\dots,1)$ point, with $L\in\Pic(X)$, $D_i=\sum_{j=1}^{m_i}c_{ij}$ effective divisors for $i=1,\dots, n-1$, $D_0=\sum_{j=1}^{m_0} c_{0j}$ where $c_{0j}\in X$ or $-c_{0j}\in X$ such that $L=\Oo(D_0)$, and $\phi_i\in H^0(\Oo(D_i))$ for $i=1,\dots, n-1$. Set (cf. \cite[\S 1]{HH} for details). 
A conjectural mirror to the $BAA$-brane $(\Ff^+,\Oo_{\Ff^+})$ over the smooth locus of the moduli space is proposed in \cite[\S6]{HH}, given by the pair  $(\M_X(n,d)^s,\Lambda_\Ff)$, where  $\Lambda_\Ff$ is a hyperholomorphic bundle supported on $\M_X(n,d)^s$ given by 
\begin{equation}
    \label{eq:hyperholo}
   \Lambda_\Ff:=\bigotimes_{i=0}^{n-1}\bigotimes_{j=1}^{m_i}\bigwedge^{n-i}\EE_{c_{ij}}.
\end{equation}
In the above, $\EE$ is the \'etale local universal bundle on $X\times \M_X(n,d)^s$, which once tensored as per \eqref{eq:hyperholo} descends to a global vector bundle \cite[\S6.2.1]{HH} on $\M_X(n,d)^s$. This bundle is endowed with a $\CC^\times$ action which makes the \'etale local projections to $X\times M_X(n,d)^s$ $\CC^\times$ equivariant and induces multiplication by $-1$ on the universal Higgs field (cf. \cite[Lemma 6.6]{HH} ). By the discussion in {\it loc. cit.},  this $\CC^\times$-equivariant structure is defined up to tensoring by a rank one $\CC^\times$- module. Now, appropriately restricting and tensoring, a  $\CC^\times$ equivariant structure is induced on $\Lambda_\Ff$, again defined up to tensorisation by a rank one module. In order to remove this ambiguity   \cite[Lemma 6.6]{HH}, one may normalise by establishing that the $\CC^\times$ equivariant structure induced on $\mathrm{det}(\bigwedge^{n-i}\EE_{c_{ij}})$ have weight
\begin{equation}
    \label{eq:normalised_weight}
w(i):=\binom{n-i}{2}-\binom{n}{2}\binom{n-1}{n-i-1}.
\end{equation}
Given another fixed point $\Ee$, mirror symmetry considerations \cite[Theorem 6.9, Remark 6.12]{HH}  lead to considering a symmetrised version of equivariant Euler pairing of the branes $(\Ee^+,\Oo_{\Ee^+})$ and $(\M_X(n,d),\Lambda_\Ff)$ (obtained by substituting $\Ss_1=\Oo_{\Ee^+}$ and $\Ss_2=\Lambda_\Ff$ in \eqref{eq:equivar_euler_pairings}) with the simpler invariant 
\begin{equation}\label{eq:equivar_euler_flows}
  m_{\Ff,\Ee}:=\chi_{\CC^\times}(\M_X(n,d),\Lambda_\Ff|_\Ee).
\end{equation}
By \cite[Remark 6.14]{HH}, if $\Ee,\Ff$ were very stable, $m_{\Ff,\Ee}$ should be a polynomial. The next result shows that, when $\Ee$ of type $(n_0,n_1)$ (and thus wobbly by Theorem \ref{thm:wobbly_comps}), $m_{\Ff,\Ee}$ is also a polynomial.
\begin{proposition}\label{prop:equiv_Euler_pairings}
        Let $\Ee\in\Fix$, $\Ff$ be as in \eqref{eq:Ff}. Then, with the above $\CC^\times$-equivariant structure on $\Lambda_\Ff$,  the pairing $m_{\Ff,\Ee}(t)$ has the following expression
    \[
    m_{\Ff,\Ee}(t)=t^{\epsilon(\Ee,\Ff)}\prod_{i<n_1}\left(\sum_{l=0}^{i}\binom{n_0}{l+n_0-i}\binom{n_1}{n_1-l}t^l\right)^{m_i}\]\[\prod_{n_1\leq i< n_0}\left(\sum_{l=0}^{n_1}\binom{n_0}{l+n_0-i}\binom{n_1}{n_1-l}t^l\right)^{m_i}\]\[\prod_{n_0\leq i\leq n-1}\left(\sum_{l=0}^{n-i}\binom{n_0}{l}\binom{n_1}{n-i-l}t^l\right)^{m_i}
    \]
    where 
    \[
    \epsilon(\Ee,\Ff):=\sum_{i=0}^{n_0-1}m_i(n_0-i)-\sum_{i=0}^{n-1}m_i((n-i)w_i),
    \]
    $w_i$ is 
    non positive and defined by the equation
    \begin{equation}
        \label{eq:wi}
        w(i)=-n_0\binom{n-1}{n-1-i}+(n-i)\binom{n}{n-i}w_i
    \end{equation}
    where $w(i)$ is defined in \eqref{eq:normalised_weight}.

    In particular $m_{\Ff,\Ee}$ is a polynomial up to a positive rational power of $t$.
\end{proposition}
\begin{proof}
   Since $(\Lambda_\Ff)|_\Ee$ is supported in dimension $0$, hence 
    \[m_{\Ff,\Ee}=\sum_{j} H^0(\M_X,\Lambda_\Ff|_{\Ee})^jt^{-j}.\]
     To compute the weight spaces, of $\Lambda_\Ff|_{\Ee}$, note that the $\CC^\times$-action is induced by one on $\EE$ yielding a weight $-1$ on the universal Higgs bundle $\Phi\in H^0(\End(\EE)\boxtimes K)$. In particular, the restriction $\EE|_{c_{ij}\times\Ee}\cong \Ee|_{c_{ij}}$ for a given fixed point $\Ee=(F_0\oplus F_1,\varphi)$ of type $(n_0,n_1)$ and $c_{ij}\in X$ satisfies that  $F_k|_{c_{ij}}$ is preserved with weight $w^k_{ij}=w^k_{i}$ is independent of $j$ and satisfies $w^0_i+w_i^1=-1$. To fully determine the weights, we note that
    \[\left(\bigwedge^{n-i}\EE_{c_{ij}}\right)|_{\Ee}=\left\{\begin{array}{ll}
    \bigoplus_{l=n_0-i}^{n_0} \bigwedge^l(F_0)_{c_{ij}}\wedge\bigwedge^{n-i-l}(F_1)_{c_{ij}}    & \textrm{ if }i<n_1; \\
     \bigoplus_{l=n_0-i}^{n-i} \bigwedge^l(F_0)_{c_{ij}}\wedge\bigwedge^{n-i-l}(F_1)_{c_{ij}}    & \textrm{ if }n_1\leq i<n_0;\\ \bigoplus_{l=0}^{n-i} \bigwedge^l(F_0)_{c_{ij}}\wedge\bigwedge^{n-i-l}(F_1)_{c_{ij}}    & \textrm{ if }n_0\leq i<n.
    \end{array}\right.\]
    Now, the dimension of $\bigwedge^l(F_0)_{c_{ij}}\wedge\bigwedge^{n-i-l}(F_1)_{c_{ij}}$ is $\binom{n_0}{l}\binom{n_1}{n-i-l}$ and there is an induced $\CC^\times$-action with weigh $lw_{i}^0+(n-i-l)w_{i}^1=-l+(n-i)w_i^1$. By multiplicativity of the Euler characteristic we find
    \[
\chi_{\CC^\times}(\M_X,\Oo_{\Ee^+},\Lambda_\Ff)=\prod_{i=0}^{n-1}\left(\chi_{\CC^\times}(\M_X,\bigwedge^{n-i}\Lambda_\Ff|_{\Ee}\right)^{m_i}=
\]
\[
\prod_{i=0}^{n-1}\left(\chi_{\CC^\times}\left(\bigoplus_{l=0}^{n-i} \bigwedge^l(F_0)_{c_{ij}}\wedge\bigwedge^{n-i-l}(F_1)_{c_{ij}}\right)\right)^{m_i}=\]\[\prod_{i=0}^{n_1-1}\left(\sum_{l=n_0-i}^{n_0} \binom{n_0}{l}\binom{n_1}{n-i-l}t^{l-(n-i)w_{i}^1}\right)^{m_i}\]\[\prod_{i=n_1}^{n_0-1}\left(\sum_{l=n_0-i}^{n-i}\binom{n_0}{l}\binom{n_1}{n-i-l}t^{l-(n-i)w_{i}^1}\right)^{m_i}\]\[\prod_{i=n_0}^{n-1}\left(\sum_{l=0}^{n-i}\binom{n_0}{l}\binom{n_1}{n-i-l}t^{l-(n-i)w_{i}^1}\right)^{m_i},
    \]
 Imposing the extra condition that the action has weight $w(i)$ on $\det(\bigwedge^{n-i}\EE_{c_{ij}})$ proves the statement about the expression of $m_{\Ee,\Ff}$ and \eqref{eq:wi}.
 
 To prove non positivity of $w_i$, it is enough to prove that $(n-i)w_i\leq0$, which is equivalent to proving non positivity of   
\begin{equation}
    \label{eq:bounding}
\binom{n}{n-i}(n-i)w_i=w(i)+n_0\binom{n-1}{n-i-1}=\binom{n-i}{2}+\binom{n-1}{n-i-1}\left(-\binom{n}{2}+n_0\right).
\end{equation}
Now, assume $n\neq n_0+1$; then $\binom{n}{2}\geq \binom{n_0+2}{2}=\frac{n_0^2+3n_0+2}{2}$, so we may bound \eqref{eq:bounding} further
\[
\leq \binom{n-i}{2}+\frac{1}{2}\binom{n-1}{n-i-1}\left(-n_0^2-n_0-2\right)\leq \binom{n-i}{2}-\frac{1}{2}\binom{n-1}{n-i-1}\left(\frac{n^2}{4}+\frac{n}{2}+2\right)
\]
where the last inequality uses that $n_0\geq n/2$. Since $w_0=0$, we may assume $i\geq 1$ to find the bound 
\[
\binom{n}{n-i}(n-i)w_i\leq \binom{n-1}{2}-\frac{1}{2}\binom{n-1}{n-2}\left(\frac{n^2}{4}+\frac{n}{2}+2\right)=\]\[=\frac{n-1}{2}\left(n-2-\frac{n^2}{4}-\frac{n}{2}-2\right)<0.
\]
We still need to check the case $n=n_0+1$, in which
\[
\binom{n}{n-i}(n-i)w_i=\binom{n-i}{2}+\binom{n-1}{n-i-1}\left(-\binom{n}{2}+n-1\right)=\]\[=\binom{n-i}{2}-\binom{n-1}{n-i-1}\binom{n-1}{2}\leq \binom{n-1}{2}\left(1-\binom{n-1}{n-i-1}\right)\leq 0.
\]
\end{proof}
\begin{remark}
    We note that $m_{\Ff,\Ee}$ is a  polynomial up to a positive rational power of $t$. This factor is the manifestation of a gerbe by which the duality should be corrected to be exact, as it happens for type $(1,\dots,1)$ pairings \cite[Theorem 7.2]{HH}. 
\end{remark}
Define
\[
m_{\Ee,\Ff}(t):=\frac{m_\Ee(t)\cdot m_{\Ff,\Ee}(t)}{m_\Ff(t)}.
\]
By \cite[Remark 6.14]{HH}, $m_{\Ee,\Ff}$ is expected to be polynomial for $\Ee, \Ff$ very stable, as in such a case  $m_{\Ee,\Ff}(t)=\chi_{\CC^\times}(\Lambda_\Ee|_\Ff)$, where $\Lambda_\Ee$ is the brane dual to  $\Oo_{\Ee^+}$. It thus provides extra obstructions to very stability of $\Ee$. 

However, in the particular case of type $(n_0,n_1)$ fixed points, an analysis of examples shows that the prediction capability for wobbliness of $m_{\Ee,\Ff}(t)$  is not better than that of the virtual equivariant multiplicities. 

\begin{example}\label{ex:euler}
    In the case $n=4$, we have for $n_0=n_1=2$
    \[
    m_{\Ee,\Ff}=[2]_t^{\delta-2(g-1)+m_1+m_2+m_3}[3]_t^{5g-5-m_2}[4]_t^{7g-7-m_1+m_2+m_3}.
    \]
    Note that $m_{\Ee,\Ff}$ is a polynomial even when $m_\Ee(t)$ isn't.
    
    For $n_0=3,n_1=1$ we find in turn
    \[
    m_{\Ee,\Ff}=[2]_t^{\delta-(g-1)+m_2+m_3}[3]_t^{5g-5-m_2}[4]_t^{7g-7-m_1+m_2+m_3}
    \]
    which is always a polynomial, just as $m_\Ee(t)$.
\end{example}

\end{document}